\documentclass[12pt]{article}

\usepackage[english]{babel}
\usepackage{amsmath}
\usepackage{lineno}
\usepackage{caption}
\usepackage{comment}
\usepackage{subcaption}
\usepackage{amssymb}
\usepackage{mathtools}
\usepackage{mathptmx} 
\usepackage{xcolor}
\usepackage{tikz}
\usepackage[square,numbers,sort]{natbib}
\usepackage{dsfont}

\usepackage{amsthm}
\usepackage{dsfont}
\usepackage[left=1in, top=1in, bottom=1in, right=1in]{geometry}
\newcommand{\derive}[2] {{\frac{\partial {#1} }{\partial {#2}}}}
\newcommand{\derd}[2] {{\frac{\vd #1 }{\vd #2}}}

\newcommand{\T}{\mathcal T}

\newcommand{\R}[0]{\mathbb{R}}

\newcommand{\abs}[1]{\left | #1\right| }

\newcommand{\norm}[1]{\left \|#1  \right \|}

\newcommand{\sset}[1]{ \{ #1\}}

\newcommand{\Leb}{\mathrm{L}}
\newcommand{\inner}{\overset{\circ}}

\newcommand{\vd}{ \mathrm{d}}
\newcommand{\intd}{\,\mathrm{d}}
\newcommand{\lo}{\mathrm{LO}}
\newcommand{\ho}{\mathrm{HO}}

\newcommand{\refs}{\mathrm{ref}}
\renewcommand{\epsilon}{\varepsilon}
\renewcommand{\phi}{\varphi}
\newcommand{\skp}[2]{\left \langle {#1}, {#2} \right \rangle}
\newcommand{\bskp}[2]{\left [ {#1}, {#2} \right ]}

\DeclareMathOperator{\argmin}{argmin}
\DeclareMathOperator{\bigO}{\mathcal{O}}
\DeclareMathOperator{\One}{\mathds{1}}
\DeclareMathOperator{\lspan}{span}
\DeclareMathOperator{\op}{\mathbb{P}}
\DeclareMathOperator{\ip}{\mathbb{I}}
\DeclareMathOperator{\ev}{\derive{U}{u}}
\newcommand{\NQ}[2]{\begin{tikzpicture}[baseline=-5.0pt, scale=0.2]
		\draw[line width=0.4mm](-1, -2) -- (-0.75, -2) node [anchor=west]{$\scriptstyle #1$} -- (-0.5, -2) -- (0.5, 2) -- (0.75, 2) node [anchor=west]{$\scriptstyle #2$} -- (1, 2);
\end{tikzpicture}}
\newcommand{\numd}{\begin{tikzpicture}[scale=0.15]
		\draw[line width=0.3mm] (1, 2.5)--(1, 0)--(0,0)--(0, 1)--(1,1);
	\end{tikzpicture}~}
\newcommand{\rskp}[2]{\left(#1\middle)\cdot\middle(#2\right)}
\theoremstyle{definition}
\newtheorem{definition}{Definition}
\theoremstyle{theorem}
\newtheorem{lemma}{Lemma}
\theoremstyle{remark}

\theoremstyle{theorem}
\newtheorem{theorem}{Theorem}
\newtheorem{remark}{Remark}

\author{Simon-Christian Klein \footnote{simon-christian.klein@tu-bs.de}}

\title{Stabilizing discontinuous Galerkin methods using Dafermos' entropy rate criterion}

\begin{document}

	\maketitle
	\abstract{A novel approach for the stabilization of the discontinuous Galerkin method based on the Dafermos entropy rate crition is presented.
		The approach is centered around the efficient solution of linear or nonlinear optimization problems in every timestep as a correction to the basic discontinuous Galerkin scheme. The thereby enforced Dafermos criterion results in improved stability compared to the basic method while retaining the order of the method in numerical experiments. Further modification of the optimization problem allows also to enforce classical entropy inequalities for the scheme. The proposed stabilization is therefore an alternative to flux-differencing, finite-volume subcells, artificial viscosity, modal filtering, and other shock capturing procedures.
}
	\section{Introduction} \label{sec:intro}
	In this work, a novel shock-capturing approach for discontinuous Galerkin (DG) schemes is proposed. The first subsection of the introduction gives a short reminder of the basic theory of hyperbolic conservation laws and the Dafermos entropy rate criterion, which is the basis for the shock-capturing technique. A second subsection covers the basic discontinuos Galerkin framework that is used as a base scheme for our construction. Section \ref{sec:theory} outlines the general idea used for the correction, while section \ref{sec:construction} and section \ref{sec:disccons} give a preliminary analysis of the properties of the modification. We close our presentation with some convincing numerical tests in section \ref{sec:NT} and an conclusion. The used notation is summarized in table \ref{tab:notation}.

\subsection{Hyperbolic conservation laws}
Systems of hyperbolic conservation laws give quantitative insight into some of the most interesting physical phenomena known to human mankind \cite{Dafermos2016Hyper}. In this work we will consider one-dimensional systems of conservation laws for $m$ conserved quantities \cite{Lax73},
\begin{equation}
	\derive{u(x, t)}{t} + \derive{f\circ u (x, t)}{x} = 0 \quad \text{for} \quad u(x, t): \R \times \R \to \R^m, \quad f: \R^m \to \R^m
	\label{eq:hpde},
\end{equation}
where $f \circ u(x, t) = f(u(x, t))$ denotes the composition of $f$ and $u$.
As classical solutions to \eqref{eq:hpde} can break down in finite time \cite{Lax73} one considers weak solutions $u$, which satisfy
\begin{equation}
	\label{eq:weaksol}
		\begin{aligned}
	\int_0^\infty \int_\R \rskp{u(x, t)}{ \derive{\phi(x, t)}{t}} &+ \rskp{ f \circ u(x, t)}{ \derive{\phi(x, t)}{x}}  \intd x \intd t \\
		&+ \int_\R \rskp{u(x, 0)}{\phi(x, 0)} \intd x = 0. 
	\end{aligned}
\end{equation}
This can be shown for classical solutions by multiplying with a suitable test function $\phi: \R \times \R \to \R^m$, integrating over the domain and using integration by parts \cite{Lax71}. Sadly, weak solutions are not unique \cite{Lax54} and one therefore hopes that additional constraints in form of entropy inequalities, 
\begin{equation}
	\derive{U \circ u}{t} + \derive{F \circ u}{x} \leq 0,
	\label{eq:eie}
\end{equation}
single out the relevant solution \cite{Lax73}. In this inequality $U: \R^m \to \R$ is the convex entropy function and $F: \R^m \to \R$ the entropy flux satisfying
\[
	\derd{U}{u} \derd{f}{u} = \derd{F}{u}.
\]
One can show that a pair of functions $(U, F)$, satisfying this relation, induces an additional conservation law for smooth classical solutions. The entropy inequality for weak vanishing viscosity solution motivates the entropy inequality, as given above \cite{Lax71}.
 Entropy inequalities are not the only criterion one can use to shrink the set of admissible solutions. Dafermos proposed an entropy rate criterion in his seminal work \cite{Dafermos72}. After defining the total entropy in the domain,
\[
E_u(t) = \int_\R U \circ u(x, t) \intd x,
\]
he conjectures that the physically relevant solution $u$ is the one satisfying 
\[
\forall t > 0: \quad \derd{E_u(t)}{t} \leq \derd{E_{\tilde u}(t)}{t}
\]
 compared to all other weak solutions $\tilde u$. See \cite{Dafermos2009MDR, Dafermos2012MD, Feireisl2014MD} for theoretical examples where the Dafermos criterion is able to reduce the amount of admissible solutions.

\subsection{Discontinuous Galerkin methods}
Discontinuous Galerkin (DG) methods, first considered for time dependent hyperbolic problems in \cite{Cockburn1989DGI, CockburnShu1989DGI}, can be derived as a generalization of classical Finite Volume (FV) methods \cite{SonarFVM}. See \cite{DG, HW2008DG, HypBookCJST} for a general introduction to DG methods. After the problem domain has been subdivided into a set $\T$ of disjoint cells one can select a basis of ansatzfunctions \[ B^T = \sset{\phi_1^T(x), \phi_2^T(x),  \dots, \phi_M^T(x)}\] associated to every cell $T \in \T$. One further conjects that the numerical solution $u(x, t)$ is a linear combination of these basis functions on every cell and only piecewise continuous in space. Discontinuities are only allowed to occur at cell interfaces. We denote the space of ansatzfunctions on cell $T$ by $V^T = \lspan B^T$. Often, the index $T$ will be omitted. One therefore writes 
\[
	u(x, t) = \sum_{T \in \T} \sum_{j \in T} u^T_j(t)\phi_j^T(x)\chi^T(x)
\]
as an ansatz for a time dependent function $u(\cdot, t)|_{T} \in V^T$ with the time dependent coefficients $u_j^T(t)$ and the characteristic functions $\chi^T(x)$ of every cell. 
A time evolution equation for our ansatz can be derived in every cell $T\ \in \T$ by multiplying the conservation law by another function $v \in V^T$, using partial integration
\begin{equation}
	\int_T \derive{v}{x} f(u^T(x, t)) \intd x - [v f(u^T(x, t))]_{T_l}^{T_r} - \int_T v\derive{u^T(x, t)}{t} \intd x = 0,
	\label{eq:dgweak}
\end{equation}
and requiring that this equation holds for all $v \in V^T$ in the same way weak solutions are defined in \eqref{eq:weaksol}. One can rewrite \eqref{eq:dgweak} as
\[
\forall T \in \T, v \in V: \skp{\derive{v}{x}}{f}_T - \bskp{v}{f}_T - \skp{v}{\derive{u^T}{t}}_T = 0
\]
using the inner products
\[
	\skp{u}{v}_T = \int_T \rskp u  v \intd x, \quad \bskp{v}{f}_T = \int_{\partial T} \rskp v f \intd O.
\]
Evaluating these inner products for $u, v \in V^T$ is done using the basis representations $(u_k)_{k=1}^M, (v_k)_{k=1}^M$ of $u$ and $v$ and the Grammian matrices $M^T$ and $S^T$ associated with cell $T$,
\[
	M_{k, l}^T = \skp{\phi_k}{\phi_l}_T, \quad S_{k, l}^T = \skp{\derive{\phi_k}{x}}{\phi_l}_T.
\]
Clearly, this is possible because of their linearity
\[
	\skp {u(\cdot, t)}{ v(\cdot, t)}  _T = \skp{\sum_k u_k \phi_k}{\sum_l v_l \phi_l}_T = \sum_{k, l}u_k \skp{\phi_k}{\phi_l}_T v_l = \rskp{u}{Mv}.
\]
A 2-norm is also induced by the inner product $\skp \cdot \cdot$ and will later play an import role
\[
	\norm{u^T(t)}_{T,2} = \sqrt{\skp{u^T(t)}{ u^T(t)}_T} = \sqrt{\rskp{u^T(t)}{M u^T(t)}},
\] 
together with the integral functional on $V^T$ 
\[
 V^T \ni u^T  \mapsto \NQ T ~ u^T(x, t) \numd x = \skp{1}{ u^T}_T = \rskp{1}{ M u^T}.
\]
Please note that the exact integration of $u^T$ depends on the inclusion of constants into $V$ \cite{GG2021RBF}, which will be assumed as given from now on.
For nonlinear fluxes the chained function $f \circ u$ is not necessarily in $V$. We therefore must find a different method to evaluate the inner products involving $f$ at least approximately. Approximating $f$ in the space $V$ via a projection $\ip_V: \R^\R \to V$ with $\ip_V = \ip_V^2$ is a common method in this case, leading to the modified scheme
\begin{equation}
	\forall v \in V: \skp{\derive{v}{x}}{\ip_V f}_T - \bskp{v}{f}_T - \skp{v}{\derive{u}{t}}_T = 0.
		\label{eq:pdg}
\end{equation}
A suitable projection is the interpolation of $f$ using suitable collocation points $\xi_k \in T$, i.e.
\[
	\ip_V f (x) = \sum_{k} \phi_k(x) f_k, \text{ with } (f_k)_{k=1}^M \in \R^M \text{ satisfying }  \forall \xi_k:  f(\xi_k) = \ip_V f(\xi_k).
\]
While one could also think of least squares projections $\op_V$, these are as hard to calculate as the sought after inner product \cite{HW2008DG}. In what follows, we use the Gauß--Lobatto points as collocation points \cite{offner2019error} and assume that the matrices $M^T, S^T$ are given with respect to the Lagrange polynomials of $\xi_k$, i.e. our basis $B$ are the Lagrange polynomials. The practical implementation should reside also to the usage of Legendre polynomials \cite{HW2008DG}.
The values of $u$ at the element boundaries are not uniquely determined and neither are their fluxes. One therefore uses monotone FV two-point fluxes as approximation of the flux over the boundary $f^* \approx f\left(\lim_{h \uparrow 0}u^T\left(x_{k+\frac 1 2} + h, t\right),  \lim_{h \downarrow 0}u^T\left(x_{k+\frac 1 2} +h, t\right)\right)$, where the limits of $u^T(x, t)$ are entered into the two arguments. As we introduced constants into our approximation space,
\[
	\skp{\One}{ \derd{u^T}{t}}_T = \underbrace{\skp{\derive{\One}{x}}{\ip_V f}_T}_{= 0} - [\One, f^*]_T = f^*_l - f^*_r,
\] holds. The DG method can therefore be seen as a FV method that uses classical FV fluxes to determine the flux of conserved quantities between cells but also advances an ansatz inside the cells using a finite element ansatz. An especially popular choice for the inter cell flux is the (local) Lax--Friedrichs flux \cite{Lax71, CockburnDG}. The complete schemes can be rewritten in vector matrix notation as
\[
	M^T \derd{u^T}{t} = S^T f(u^T(t)) - \begin{pmatrix} \phi^T_1(x_r)f^*_r-\phi^T_1(x_l) f^*_l \\ \vdots \\  \phi^T_N(x_r) f^*_r  - \phi^T_N(x_l) f^*_r\end{pmatrix}.
\]
The resulting scheme can be integrated int time, for instance, by strong stability preserving Runge--Kutta methods \cite{SO1988, SSPRK} or any other solver for ordinary differential equations. We will not consider the problem that if the semidiscrete scheme satisfies some bound or property will the numerical solution to the ODE system in general \emph{not} satisfy this bound. Refer to \cite{ranocha2020strong} and \cite{OFFNER202015} for an outline of the problems a set of possible solutions.
Classical DG type methods use total variation limiters \cite{CockburnDG} and troubled cell indicators that switch a cell into a first order FV mode \cite{CockburnDG} or reconstruct the ansatz from neighboring cells to enforce stability and improve their abilities to calculate discontinuous solutions. More recent methods using the summation by parts (SBP) property \cite{Gassner, ranocha2016summation, ShuDGReview} in conjunction with flux differencing \cite{Fisher2013Fluxdiff} allow to obtain stability results without the usage of limiters. Another option  are the techniques of Abgrall, who introduced a general entropy correction term \cite{abgrall2019reinterpretation}. While the aforementioned results provide entropy stability, the resulting schemes often need additional stabilization to calculate solutions containing shocks. We will design two similar entropy correction techniques in this work for the usage in the DG framework based on Dafermos' entropy rate criterion that also serves as a shock-capturing technique. Prior usage of the Dafermos entropy rate criterion was based on solving the variational problem given by the entropy rate criterion in a convex subset of the fluxes \cite{klein2022using}. In this work, the entropy rate criterion will be applied to the approximate solutions inside every cell of the DG scheme while the fluxes between cells will be left untouched. Therefore existing reasonable numerical fluxes can still be used. One should further note that the presented correction operations can be generalized to multiblock SBP schemes and the author will consider these methods in the future. It should be noted that the presented techniques will be tested here only without SBP-SAT and flux-differencing techniques to underline their abilities. 

\begin{table}
\begin{tabular}{c | c}
	A cell of the subdivision $\T$ of the domain $\Omega$ & $T$ \\
	The left and right boundaries of cell $T$ & $T_l$, $T_r$\\
	The space of ansatzfunctions for an cell $T$& $V^T$ \\
	$\Leb^2$ projection of $u$ onto $V$ & $\op_V u$.\\
	Interpolation of $u$ on $V$ w.r.t. the collocation points $(\xi_k)_{k=1}^N$& $\ip_V u$ \\ 
	Vector of nodal values in cell $T$ at time $t$& $u^T(t)$ \\
	Ansatz function in cell $T$ at position $x$ and time $t$ & $u^T(x, t)$ \\
	The total entropy in cell $T$ & $E_{u, T}(t)$ \\
	The discrete total entropy in cell $T$ & $E^T(t)$ \\
	Entropy variables in cell $T$ & $\derive{U}{u}(u^T(x, t))$ \\
	Vector of nodal values of the entropy variables in cell $T$ at $t$ & $\ev^T(t)$ \\
	Interpolation of the entropy variables in cell $T$ on $V$& $\ev^T(x, t)$ \\
	Numerical quadrature of $f(x)$ on the intervall $[a, b]$ & $\NQ{a}{b} f(x) \numd x = \sum_k \omega_k f(x_k)$ \\
	The canonical inner product between $a, b \in \R^n$ & $\rskp{a}{b}$ \\
	The inner product between $u$ and $v$ on cell $T$ & $\skp{u}{v}_T$ \\
	The $p$ norm of $u$ in cell $T$ & $\norm{u}_{T, p}$ \\
	Exact solution to the initial condition $u(x, t_0)$ after $t-t_0$ & $H(u(\cdot, t_0), t-t_0)$. \\
	Mean value of subcell $k$ of $N$, $u^T$ as initial condition & $u^{T, N}_k$ \\
	
\end{tabular}
\caption{Used notation. As a general rule, quantities with only $t$ as an argument are a vector of nodal values at a certain time. Values with $x$ and $t$ in their argument list are functions that were evaluated at these values. Objects with $T$ added as exponent are approximations of the quantity in the cell $T$.}
\label{tab:notation}
\end{table}
	\section{Dafermos modification of DG methods} \label{sec:theory}
	In the last chapter the DG schemes were introduced. Let us denote an assumed exact solution operator as $H(u(\cdot, t_0), \tau)$ that maps an initial condition $u(\cdot, t_0)$ onto the solution at time $t_0 + \tau$. The finite dimensional vector space $V$ extorts several approximations, like the projection $\ip_V$ of nonlinear fluxes and that the time derivative of the coefficients $\derd{u^T}{t} \in V$ satisfies the weak formulation only for $v \in V$. The time discretization is therefore in general different from the projection of $\derive{H(u(\cdot, t), \tau)}{\tau}$ onto $V$ and also this projection differs from $\derive{H(u(\cdot, t), \tau)}{\tau}$.
We therefore select a trajectory from the initial condition leading away from the exact solution $u$ or its projection onto $V$. One further knows that weak solutions, which are standard for hyperbolic conservation laws and built into the DG method, allow more than one path for an exact weak solution. We will therefore modify the outlined vanilla DG scheme to rule out unwanted deviations from the assumed trajectory of the assumed exact entropy solution. Unwanted deviations are deviations that lead to
instabilities,
non-admissible solutions,
oscillations or
high approximation errors.
It is therefore conjectured that changing the trajectory $\derd{u^T(t)}{t}$ of our approximate solution to satisfy a modified Dafermos entropy rate criterion performs the needed deviation to rule out these problems.
We begin the construction of our modification by the definition of a per cell total entropy as
\[
	E_{u, T}(t) = \int_T U(u^T(x, t)) \intd x = \skp{\One}{ U(u^T(x, t))}_T \approx \skp{\One}{ U^T(t)}_T = E^T(t),
\]
where $U^T(t) = U(u^T(t))$ is the vector composed of the values of $U(u)$ evaluated with the vector of nodal values $u^T(t)$. 
 We will see that our numerical entropy functional for DG schemes is (strictly) convex and localy Lipschitz continuous under certain circumstances. One of the needed ingredients is a quadrature formula with perfect stability \cite{Glaubitz2021Quad, Glaubitz2021Cub}.
\begin{definition}[Perfectly Stable Cubature \cite{Glaubitz2021Quad, Glaubitz2021Cub}]
	A quadrature or cubature formula $\omega_k$ is termed perfectly stable if 
	\[
		\forall k: \omega_k \geq 0
	\]
	holds and the cubature formula is exact for constants, i.e. \( \sum_k \omega_k = \mu(\Omega)\), where $\mu(\Omega) = \int_{\Omega} 1 \intd x$ shall be the measure of the integrated region $\Omega$.
\end{definition}
We will now state some useful properties of our discrete entropy functional condensed into the following lemma. 
\begin{lemma}
	Let $\omega_k$ be a perfectly stable quadrature formula on the cell $T$. Then the nonlinear functional
	\[
		E^T(v) = \sum_{k} \omega_k U(v(x_k))
	\] is local Lipschitz continuous. If $U$ is a strictly convex entropy, then $E^T$ is strictly convex. Further, the entropy functional is exact for constant functions and an upper bound for the entropy of the mean value in the cell
	\[
		U(\overline u) = U\left( \frac{\int u_T \intd x}{\mu(T)} \right) \leq \frac{E^T(u)}{\mu(T)},
	\] as used as per cell entropy in FV methods.
	
	\begin{proof}
		We begin by showing that the entropy functional is convex, if $U$ is convex. Let $\lambda \in [0,1]$ and $u, v \in V$ be arbitrary. The perfect stability of the quadrature implies
		\[
		\begin{aligned}
					E^T(\lambda u + (1-\lambda) v) &= \sum_k \omega_k U(\lambda u(x_k) + (1-\lambda) v(x_k)) \\
					&\leq \lambda \sum_k \omega_k U(u(x_k)) + (1-\lambda) \sum_k \omega_k U(v(x_k)) \\
					&= \lambda E^T(u) + (1-\lambda) E^T(v).
				\end{aligned}
		\]
		If $U$ is even strictly convex so is also $E^T$ because in this case »strictly less« can be put in place of »less« in the above derivation if $u \neq v$.
		One further knows that $E^T$ is local Lipschitz continuous in the $\norm{\cdot}_\infty$ norm with constant $L = L_U$ if $U$ is local L-Lipschitzcontinuous with constant $\mu(T)L_U$
		\[
		\begin{aligned}
			\abs{E^T(u) - E^T(v)} \leq& \sum \omega_k \abs {U(u(x_k)) - U(v(x_k)} \\
			\leq& \sum \omega_k L_U\norm{u-v}_\infty = \mu(T)L_U\norm{u-v}_\infty .
			\end{aligned}
		\]
		To show that the functional is exact for constants we just note that our quadrature is exact for constants and therefore
		\[
			E^T(\bar v) = \sum_k \omega_k U(\bar v(x_k)) = \sum_k \omega_k U(\bar v) = U(\bar v)
		\]
		holds. The last property follows because the quadrature is exact for constants and all $v \in V$ and the calculation of the mean value using this quadrature results in a convex combination of point evaluations
		\[
			U(\bar v) = U\left(\frac{\sum_k \omega_k v(x_k)}{\sum_k \omega_k} \right) \leq \frac{\sum_k \omega_k U(v(x_k))}{\sum_k \omega_k} = \frac{E^T(v)}{\mu(T)}.  
		\]
		\end{proof}
\end{lemma}

Clearly, one now aims to deviate into the reduction of this nonlinear but convex functional, or semidiscrecetly stated to minimize the derivative of this functional with respect to time
\[
	\derd{E_{u, T}(t)}{t} = \int_T \derd{U}{u} \derive{u}{t} = \skp{\derd U u}{\derive u t}_T \approx \skp {\derd {U^T} u}{\derd{u^T}{t}}_T \approx \derd{E^T}{t},
\]
towards the smallest values by varying $\derd{u^T}{t}$. It should be stressed, that our numerical approximation of the time derivative of $E_{u, T}$ is \emph{not equivalent} to the time derivative of our numerical approximation of $E_{u, T}$, as the mass matrix $M$ is in general not diagonal
\[
	\skp{\derd {U^T} u}{ \derd {u^T} t}_T = \skp{\derd {U^T} u}{M \derd {u^T} t} \neq \skp{\One}{M\left(\derd {U^T} u \derd {u^T} t\right)} = \skp{\One}{\derd {U^T} u \derd {u^T} t}_T.
\]
But as we are interested in discretizing the steepest entropy descent, and not the steepest discrete entropy descent, the first form will be used.
Our modification will result in the following scheme in the discrete case in every step.
\begin{enumerate}
	\item Given a state $u^{T}(t^n)$ calculate the next step $\tilde u^{T}(t^{n+1})$ using a vanilla RK-DG method as outlined in the last chapter.
	\item Calculate an error estimate $\delta^T$ for the solution $\tilde u^{T}(t^{n+1})$, i.e. an reasonably small $\delta^T \geq 0$ with the property
	\[
		\norm{H(u^T(t^n), \Delta t) - \tilde u^T(t^{n+1})}_T \leq \delta^T,
	\]
	\emph{without} complete knowledge of the exact solution operator $H(\cdot, \Delta t)$.
	\item Solve the optimization problem
	\begin{equation}
		\begin{aligned}
		u^{T}(t^{n+1}) &= \argmin_{u \in Z} E^T(u) \text{, with } \\
			Z &= \left\{u \in \R^n \middle |  \langle 1, u \rangle_T = \langle 1, \tilde u^T(t^{n+1}) \rangle_T \wedge\norm{u-\tilde u^T(t^{n+1})}_T \leq \epsilon(\delta^T) \right \},
		\end{aligned}
		\label{eq:discprob}
	\end{equation}
		where $\epsilon$ is a given function of the error $\delta$ and $u^T(t^n)$.
\end{enumerate}
This is a direct statement of the modified Dafermos theorem pioneered by the author in \cite{klein2022using}. The original Dafermos entropy rate criterion is thereby augmented with the two additional constraints:
\begin{itemize}
	\item The introduced additional approximation error by the dissipation should be small, i.e. the dissipation should relate to the assumed error towards the exact solution. We would like to enforce $\lim_{\delta \to 0} \epsilon = 0$.
	\item The resulting discretisation should still be conservative as defined in \cite{shi2018local}.
	\end{itemize}
These two properties are engraved into the statement that the solution of the optimization problem should still have the same mean value in the cell, i.e. the basic FV method is unchanged, and that the error introduced by the dissipation is smaller than the bound $\epsilon(\delta)$. We will therefore call this scheme Dafermos Runge--Kutta Discontinuous Galerkin (DRKDG). While the exact solution of this problem in every time step is not feasible it should be noted that, if a strictly convex entropy $U$ is used, the resulting optimization problem has an unique solution,  c.f.~section \ref{sec:disccons}. A related algorithm also can be given for a semi-discrete scheme in the following form.
\begin{enumerate}
	\item Calculate the time derivative $\derd{\tilde u^T}{t}$ of the semidiscretisation using a vanilla DG scheme.
	\item Estimate the error $\delta^T$ in cell $T$, i.e. find $\delta^T \geq 0$ with
	\[
		\norm{\derd{\tilde u^T}{t} + \derive{f(u^T(x, t))}{x}}_T \leq \delta^T
	\]
	in a suitable norm. In this case $\derive{f(u^T(x, t))}{x}$ must not be understood as a classical derivative, as it will exist only in a distributional sense.
	\item Solve the optimization problem 
		\begin{equation}
			\label{prob:semi}
			\begin{aligned}
			\derd{u^T}{t} &= \argmin_{u_t \in Z} \derd{E^T}{t}\left(u^T, u_t \right) \text{ with } \\
			 Z &=\left\{u_t \in \R \middle| \skp{ 1}{ u_t}_T = 0  \wedge \norm{u_t - \derd{\tilde u_{k}}{t}}_T \leq \epsilon(\delta^T) \right \},
			\end{aligned}
		\end{equation}
		where $\epsilon(\delta^T)$ is a given function of the error $ \delta^T$ and $u^T(t)$.
\end{enumerate}
We will call this second version Dafermos Discontinouous Galerkin (DDG). The optimization problem in this modified algorithm can be solved exactly, as we will see in the next chapter.
Clearly, finding the correct error estimate and solving the optimization problems are intricate steps in the algorithm and we will often reside to using the semi-discrete form of the algorithm because the optimization problem, and the error estimate, are significantly simpler. Another motivation for this method can be made by the following observation. Assume, that the time derivative of the exact solution $\derive{u}{t} = \derive{H(u^T(\cdot, t), \tau)}{\tau}$ and the entropy variable $\derd{U(u)}{u}(u^T(\cdot, t))$ lie in the approximation space $V$. This implies $\derive{u}{t}$ exists almost everywhere. Then follows
\[
	\int_T \derive{U(u^T(x, t))}{t} \intd x =  \skp{\derd{U}{u}}{\derive{ u}{t}}_T \leq F^T_l - F^T_r 
\]
in the sense of distributions in time for an entropy stable numerical FV flux $f$ with numerical entropy flux $F$, from an integration of the entropy inequality \eqref{eq:eie} and the chain rule at all points where $\derive{u}{t}$ exists. Assume now that $\derd{u^T}{t}$ is the time derivative of a (numerical) approximation of $u(x, t)$ inside $T$ that satisfies 
\[
	\norm{\derive u t - \derd{u^T}{t}}_{T, 2} \leq \delta.
\]
The Cauchy--Schwarz inequality \cite{LaxFun} allows us to bound the entropy production in this case
\[
\begin{aligned}
	\skp{\derd{U}{u}}{\derd{u^T}{t} - {\derive{u}{t}}}_T \leq \norm{\derd{U}{u}}_{T, 2} \delta 
	\implies \skp{\derd{U}{u}}{\derd{u^T}{t}}_T \leq \skp{\derd{U}{u}}{{\derive{u}{t}}}_T + \delta \norm{\derd{U}{u}}_{T, 2}
	\end{aligned}
\]
against the entropy production of the exact solution. If we now apply an entropy correction to $\derd{u^T}{t}$ by adding a deviation of length $\delta$ into the direction of the steepest entropy descent follows
\[
	\begin{aligned}
	\skp{\derd{U}{u}}{\derd{u^T}{t}  - \delta \frac{\derd{U}{u}}{\norm{\derd{U}{u}}}}_{T}  &= \skp{\derd{U}{u}}{\derd{u^T}{t}}_T - \delta \norm{\derd{U}{u}}_{T, 2}\\
	 &\leq\skp{\derd{U}{u}}{{\derive{u}{t}}}_{T} + \delta \norm{\derd{U}{u}}_{T, 2}- \delta \norm{\derd{U}{u}}_{T, 2} \\
	  &\leq F^T_l - F^T_r,
	\end{aligned}
\]
that the entropy production of our numerical solution is bounded by the entropy production of the exact solution. This would encourage us to use $\epsilon(\delta) = \delta$ in the previous algorithm.
Sadly, in general, we will not have an exact integration of the inner products as the exact solution and entropy will not lie in our space $V$. We are further restricted to entropy descent directions that do not change the mean value of the cell. The existence and regularity of $u$ for multidimensional systems is part of ongoing research and a generalization of this observation to counteract these problems will be part of the next section. 
	\section{Construction of the semidiscrete scheme} \label{sec:construction}
We begin our search for error estimators with the design of an error estimator for the semi-discrete scheme, i.e., the error of the derivative $\derd{u^T}{t}$ in relation to $\derive{u}{t}$. We will afterwards generalize this semi-discrete case to design an error estimator for the discrete scheme.
Estimating the error in FV methods was already done using linearized conservation laws in the sense of Friedrichs systems \cite{Sonar1993SWNorm, SA1998PEI,SS1998DGN,IS2001SORE} and numerical integration of the residual. We will use a similar procedure as error estimate for our semi-discrete problem, but will directly apply a quadrature to the error between our weak solution and the projection of a second suitably fine solution onto the space $V$ of ansatz functions. This seems problematic first, and two questions arise. The first one is why to project a function if we could just use this function as the solution of our scheme, and the second concerns the cost of calculating this more accurate solution. Our more reference will be the result of a subcell scheme \cite{MS2014FVS}. We will see that the projection of this subcell scheme depends on the size $h$ of the subcells in such a way that one easily goes over to the limit $h \to 0$. This limit will be a lot cheaper than the subcell scheme it replaces. We will first define our subcell scheme in the following way
\begin{definition}[Finite Volume subcells]
	Let $u^T(x)$ be an ansatz on a cell $ T = [x_l, x_r] \subset \R$ and let \[
		x_l = x_{\frac 1 2} < x_1 < x_{\frac 3 2} < \dots < x_{N} < x_{N+ \frac 1 2} = x_r\] be a subdivision of this cell into $N$ subcells around $x_k$ and with left boundary $x_{k-\frac 1 2}$ and right boundary $x_{k+\frac 1 2}$. Let us denote the space of piecewise constant functions on $T$ with discontinuities at $x_{k+\frac 1 2}$ by
       \[
			P^0_{(x_k)_k} = \{f: \R \to \R | \forall k \in \{1, \ldots, N \}, \forall x \in [x_{k-\frac 1 2}, x_{k+\frac 1 2}]: \, f(x) = f_k  \}.
		\]
		As a projection of our ansatz into this space shall the calculation of mean values
		\[
			u^{T, N}_k = \frac{1}{x_{k+\frac 1 2} - x_{k - \frac 1 2}} \int_{x_{k-\frac 1 2}}^{x_{k+\frac 1 2}} u^T(x, t) \intd x
		\]
		be used.
		We define by 
		\[
			G: V \to P^0_{(x_k)_k}, u^T \mapsto \sum_{k=1}^N \chi_{[x_{k-\frac 1 2}, x_{k+\frac 1 2}]} u^{T, N}_k
		\]
		the projection of our ansatz function onto the space of piecewise constant functions.
\end{definition}
If an arbitrary $u^T$ is projected in this way onto a function in  $P^0_{(x_k)_k}$ it can be used as an initial condition to a low order FV scheme on the subcells. The used numerical flux can be the same as the one used as intercell flux of the DG scheme. The solution of this scheme will, at least for small times, be a function in $\Leb^p(\Omega)$ at any fixed instant of time. We can moreover interpret it as a differentiable mapping
\[
	u^{T, N}(\cdot, t): \R \to \Leb^p, t \mapsto \sum_{k=1}^N \chi_{[x_{k-\frac 1 2}, x_{k+\frac 1 2}]}(x) u^T_k(t)
\]
 by the following definition. The result at any fixed time lends itself to be projected back onto the space $V$.
\begin{definition}[Semi-discrete schemes as $\Leb^p$ functions]
	Let $f(u_l, u_r)$ be an entropy stable two point flux. We will interpret to a given initial state $u^T_k$ the solution $u^{T, N}_k(t)$ of the semi-discrete scheme
	\[
		\derd{u^{T, N}_k}{t} = \frac{f\left(u^{T, N}_{k-1}, u^{T, N}_{k}\right) - f\left(u^{T, N}_{k}, u^{T, N}_{k+1}\right)}{x_{k+\frac 1 2} - x_{k-\frac 1 2}}
	\]
	as a function $u^{T, N}(\cdot, t) \in \Leb^p(\Omega)$ by defining
	\[
		u^{T, N}(x, t) = \sum_k  \chi_{\left[x_{k-\frac 1 2}, x_{k+\frac 1 2}\right]}(x) u^{T, N}_k(t).
	\]
	We therefore also can give an interpretation of the time derivative of the scheme as an $\Leb^p$ function
	\[
		\begin{aligned}
		\derive{u^{T, N}(x, t)}{t} 
		&= \sum_k  \chi_{\left[x_{k-\frac 1 2}, x_{k+\frac 1 2}\right]} \frac{f\left(u^{T, N}_{k-1}, u^{T, N}_{k}\right) - f\left(u^{T, N}_{k}, u^{T, N}_{k+1}\right)}{x_{k+\frac 1 2} - x_{k+\frac 1 2}} \\
		&= \sum_k \chi_{\left[x_{k-\frac 1 2}, x_{k+\frac 1 2}\right]} \derd{u^{T, N}_k(t)}{t}.
		\end{aligned}
	\]
\end{definition}
	The last definition allows us to devise approximate solutions to a given ansatz $u^T$ as initial state by choosing an appropriate numerical flux $f(u_l, u_r)$, and we will often do so by choosing the same flux as in the DG scheme under consideration, projecting, solving, and projecting back. The limiting process $N \to \infty$ is delicate if the initial condition is discontinuous between cells, even at $t = 0$. The $\Leb^1$ norm of the derivative for the cell stays bounded while all norms with higher exponents blow up. This can be seen for the $\Leb^2$ norm for example as
	\[
		\begin{aligned}
		\norm{\derive{u^{T, N}(t)}{t}}_{T, 2} =& \sqrt{\sum_{k=1}^N \left(\frac{f\left(u^{T, N}_{k-1}, u^{T, N}_{k}\right) - f\left(u^{T, N}_{k}, u^{T, N}_{k+1}\right)}{x_{k+\frac 1 2} - x_{k+\frac 1 2}}\right)^2 (x_{k+\frac 1 2} - x_{k-\frac 1 2})} \\
		\geq& \frac{\abs{f\left(u^{T, N}_{0}, u^{T, N}_{1}\right) - f\left(u^{T, N}_{1}, u^{T, N}_{2}\right)}}{\sqrt{x_{1+\frac 1 2} - x_{\frac 1 2}}} \xrightarrow{N \to \infty} \infty
		\end{aligned}
	\]
	holds as a consequence of $x_{k + \frac 1 2} \to x_{k - \frac 1 2}$ in this case, because the fluxes are discontinuous at the edges. For the inner of the domain and also the edges, if $u^T(x, t)$ is differentiable at the cell interfaces, follows on the contrary by the consistency of the flux
	\[
		\frac{f\left(u^{T, N}_{k-1}, u^{T, N}_{k}\right) - f\left(u^{T, N}_{k}, u^{T, N}_{k+1}\right)}{x_{k+\frac 1 2} - x_{k+\frac 1 2}} \xrightarrow{N \to \infty} -\derive{f(u^T)}{x}\left(\lim_{N \to \infty} x_k\right).
	\]
	We will split the derivative of our approximate solution from now on into two parts which we will call the singular and the regular part. The singular part
	\[
	\begin{aligned}
		(1-R)\derive{u^{T, N}}{t} 
		=&\frac{f\left(u^{T, N}_{0}, u^{T, N}_{1}\right) - f\left(u^{T, N}_{1}, u^{T, N}_{2}\right)}{x_{\frac 3 2} - x_{\frac 1 2}} \chi_{\left[x_{\frac 1 2}, x_{\frac 3 2}\right]}(x) \\
		 &+ \frac{f\left(u^{T, N}_{N-1}, u^{T, N}_{N}\right) - f\left(u^{T, N}_{N}, u^{T, N}_{N+1}\right)}{x_{N+\frac 1 2} - x_{N-\frac 1 2}} \chi_{\left[x_{N-\frac 1 2}, x_{N + \frac 1 2}\right]}(x)
	\end{aligned}
	\] is the derivative of the subcells next to the cell boundary of the outer big cell. The regular part
	\[
		R \derive{u^{T, N}}{t} = \sum_{k=2}^{N-1} \chi_{\left[x_{k-\frac 1 2}, x_{k+\frac 1 2}\right]}\frac{f\left(u^{T, N}_{k-1}, u^{T, N}_{k}\right) - f\left(u^{T, N}_{k}, u^{T, N}_{k+1}\right)}{x_{k+\frac 1 2} - x_{k-\frac 1 2}} 
	\]
	 shall be the derivative of the subcells using only the extrapolated inner polynomial of the cell. Although the $2$-norm of our approximate solution derivative blows up under grid refinement, we can still project this approximate solution onto the function space of our ansatzes $V$ for fixed $N$. One can even find a closed expression for the limit of this projection for $N \to \infty$, and this limit exists as the sequence of projections is bounded in the $\Leb^2$ norm by the norm of the sequence of solutions prior to the projection in the $\Leb^1$ norm
	\[
		\norm{\op_V \derive{u^{T, N}(x, t)}{t}}_{T, 2} 
		\leq
		\sum_{k = 1}^M \norm{\frac{\phi_k}{\skp{\phi_k}{\phi_k}}_{T, 2} \skp{\phi_k}{\derive{u^{T,N}}{t}}}_{T, 2} 
		\leq 
		\sum_{k = 1}^M \frac{\norm{\phi_k}_{T, \infty} \norm{\derive{u^{T, N}}{t}}_{T,1}}{\norm {\phi_k}_{T,2}}.
	\]
	Because of our split into the regular and singular part we can estimate the norm of the time derivative as
	\[
		\norm{\derive{u^{T, N}}{t}}_{T, 1} \leq \norm{R\derive{u^{T, N}}{t}}_{T, 1} + \norm{(1-R)\derive{u^{T, N}}{t}}_{T, 1}.
	\]
	The norm of the regular part has to converge to $\derive{f}{x}$, as stated before, and therefore be bounded. The singular part is clearly bounded by 
	\[
		\begin{aligned}
		&\norm{(1-R)\derive{u^{T, N}}{t}}_{T, 1} \\
		&\leq \abs{f\left(u^{T, N}_{0}, u^{T, N}_{1}\right) - f\left(u^{T, N}_{1}, u^{T, N}_{2}\right)}+ \abs{f\left(u^{T, N}_{N-1}, u^{T, N}_{N}\right)- f\left(u^{T, N}_{N}, u^{T, N}_{N+1}\right)}.
		\end{aligned}
	\]
	
	We will now calculate a closed expression for the limit $N \to \infty$ for the sequence of these projections. For the regular part  
	\[
		R \derive{u^{T, N}}{t} \xrightarrow{N \to \infty} -\derive{f}{x} \implies \op_V R \derive{u^{T, N}}{t} \xrightarrow{N \to \infty} - \op_V \derive{f}{x}
	\]
	follows because of the convergence in the $\Leb^2$ norm of the regular part towards the strong spacial flux derivative. Calculating the $\Leb^2$ projection of the singular part boils down to the calculation 
	\[
		\begin{aligned}
		&&\int_T \phi_k(x) \chi_{\left[x_{\frac 1 2}, x_{\frac 3 2}\right]} \frac{f\left(u^{T, N}_0, u^{T, N}_1\right) - f\left(u^{T, N}_{1}, u^{T, N}_{2}\right)}{x_{\frac 3 2} - x_{\frac 1 2}}  \intd x \\
		&=& 
		\int_{x_{\frac 1 2}}^{x_{\frac 3 2}} \phi_k(x) \frac{f\left(u^{T, N}_{0}, u^{T, N}_{1}\right) - f\left(u^{T, N}_{1}, u^{T, N}_{2}\right)}{x_{\frac 3 2} - x_{\frac 1 2}} \intd x \\
		&=& 
		\phi_k(\xi) \left(f\left(u^{T, N}_{0}, u^{T, N}_{1}\right) - f\left(u^{T, N}_{1}, u^{T, N}_{2}\right)\right), \quad \xi \in \left[x_{\frac 1 2}, x_{\frac 3 2}\right] \\
		&&\xrightarrow{N \to \infty} \phi_k(0) \left(f^* - f\left(u^T\left(x_l\right)\right)\right). 
		 \end{aligned}
		 \]
		 This shows
		 \[
		  \op_V(1-R) \derive{u^{T, N}}{t} \xrightarrow{N \to \infty} \sum_k \frac{\phi_k(x_l)}{\skp{\phi_k}{\phi_k}} \left(f^*_l - f\left(u^T\left(x_l\right)\right)\right) + \frac{\phi_k(x_r)}{\skp{\phi_k}{\phi_k}} \left(f\left(u^T\left(x_r\right)\right) - f^*_r\right).
	\]
	Here we used the mean value theorem of integration.
	Interestingly, the sequence of the projected approximate solution derivatives does not blow up in the norm but stays bounded also in all norms as our ansatz functions are bounded in all norms. We will use this knowledge later to calculate the distance between this projected approximate solution and the solution yielded by the DG method as an error estimate.
	 This is possible as we can calculate the spatial derivative of our ansatz function $u^T(x, t)$, that is smooth in the inner of every cell, exactly as
\[
	\derive{u^T(x, t)}{x} = \derive{~}{x} \sum_{l=1} u^T_l(t)\phi_l(x) = \sum_{l=1} u^T_l(t) \derive{\phi_l(x)}{x}.
\]
One has therefore for $x \in \inner T$
\[
	\derive{u}{t} = - \derive{f(u^T(x, t))}{x}= - \derive{f}{u}\derive{u^T}{x}
\]
by usage of the chain rule and has therefore, in theory, access to the projection of the limit of the low order scheme.

Sadly, the exact evaluation of the error is still a nontrivial task as we would need to project the regular part of our assumed exact solutions derivative, calculated using the strong form, onto our vectorspace $V$. Because this is in general a function $u \in C^k$ is this a hard problem. Incidentally, the singular part can be projected easily, as we saw earlier. We will therefore estimate the norm of the distance between the limit of the projections and the DG scheme by
\[
\begin{aligned}
 &\norm{\derd {u^T} t - \op_V R \derive{u^{T, N}} t - \op_V (1-R) \derive {u^{T, N}} t}_p \\
= &\norm{\op_V\left(\derd {u^T} t - R \derive{u^{T, N}} t - \op_V(1-R) \derive {u^{T, N}} t \right)}_p \\
\leq& \norm{\derd {u^T} t - R \derive{u^{T, N}} t - \op_V (1-R) \derive {u^{T, N}} t}_p,
\end{aligned}
\]
where the estimate follows from the fact that $\derd {u^T} t$ and $\op_V (1-R) \derive{u^{T, N}} t$ lie in $V$, and are therefore fixed points of the projection, and the projection operator has a norm smaller than one.    The author would like to stress that the function in the norm on the right hand side of the equation is as smooth as the numerical flux function and can be evaluated exactly at any point. It is therefore only logical to evaluate this norm in case of the $2$ norm using numerical quadrature, which was done in the numerical tests. 
 Quadrature rules of the Gauß--Legendre family were chosen in the numerical tests. We will from now on refer to the limit of the regular part of the approximate solution by the low-order scheme together with the limit of the irregular part projected onto $V$ as the reference solution
\[
	\derive{u^\refs}{t} = R \derive{u^{T, \infty}} t - \op_V (1-R) \derive {u^{T, \infty}} t.
\]

As explained before, our solver needs to solve an optimization problem after every timestep for the discrete, or for every evaulation of the time derivative in the semidiscrete algorithm. We will first consider the optimization problem of the semidiscrete scheme as this will in fact also be a building block for our approximate solution to the optimization problem in the discrete case. The semidiscrete problem is, luckily, a linear one and can be solved exactly, as we will see in the next lemma.

\begin{lemma}[Restricted gradien descent] \label{lem:rgd}
	A solution $s$ of the semidiscrete optimization problem \eqref{prob:semi} is given by
	\[
		g = -\derd{U^T}{u}(u^T), \quad h = g - \frac{\skp {\One}{g}_T}{\skp{\One}{\One}_T} \One, \quad s =\begin{cases} \frac{\epsilon h}{\sqrt{\skp{h}{h}_T}} & h \neq 0 \\ 0 & h = 0\end{cases}
	\]
	and this solution is unique if $\derd {U^T} u \neq c \One$ holds, i.e. when the solution in the cell is not constant.
	\begin{proof}
	We beginn our proof by showing that $s \in Z$ holds. Clearly,
	\[
		\skp{\One}{ s}_T = \frac{\epsilon}{\sqrt{\skp{h}{h}_T}} \left(\skp{\One}{g}_T - \skp{\One}{g}_T\frac{\skp{\One}{\One}_T}{\skp{\One}{\One}_T} \right) = 0
	\]		
	shows that the solution lies in the linear subspace $V \supset W = \{v \in V|\skp{\One}{v} _T= 0\}$. The rescaling from $h$ to $s$ also implies
	\[
		\norm{s}_{T, 2} = \sqrt{\skp{s}{s}_T} = \epsilon \frac{\sqrt{\skp{h}{h}_T}}{\sqrt{\skp{h}{h}_T}} = \epsilon.
	\]
	The set $Z = W \cap B_\epsilon$ is clearly the intersection of vectors of length less than or equal $\epsilon$ and the subspace $W$ and hence $s$ lies inside this vectorspace.
	After the admissibility of $s$ was established we can take care of the optimality. We will first show that for the restriction $g \in B_\epsilon$ is a scaled version of $g = -\derd{U^T}{u}$ indeed the optimal descent direction. We will afterwards show that the optimal descent direction is $s = \op_w g$ if only directions in the subspace $W$ are considered. Clearly,
	\[
		\skp{\derd {U^T} u}{\epsilon \frac{g}{\norm {g}_{T, 2}}}_T = \skp{\derd {U^T} u}{-\epsilon \frac{\derd {U^T} u}{\norm{\derd {U^T} u}_{T, 2}}}_T = - \epsilon \norm{\derd {U^T} u}_{T,2}
	\] 
	holds, because $\derd {U^T} u$ and $s$ are colinear.
	Let now $v \in V$ be arbitrary with $\norm{v}_{T, 2} = \epsilon$. The Cauchy-Schwarz inequality implies
	\[
		\abs{\skp{\derd{U^T}{u}}{v}_T} \leq \norm{\derd  {U^T} u }_T \epsilon.
	\]
	It therefore follows
	\[
		\skp{\derd {U^T} u}{s}_T \leq \skp{\derd {U^T} u}{v}_T, 
	\]
	and this inequality is strict for $\derd {U^T} u \neq c \One$ and $v \neq s$, as in this case $s \neq 0$ and $\derd {U^T} u \neq0$ follows and the Cauchy-Schwarz inequality is equal only if $\derd {U^T} u$ and $s$ are colinear in this case.
	 We note in passing that $\skp{\cdot}{ \cdot}_T$, as a inner product on $V$, also is an inner product on $W$ and hence the Cauchy-Schwartz inequality applies for all elements in $W$. We can further decompose $V$ into $W^\perp$ and $W$ and if at least one vector $v \in W$ of two vectors $v, w \in V$ is from $W$ follows
	\[
		\skp{v}{w} = \skp{\op_W v + \op_{W^\perp}v}{w} = \skp{\op_W v}{w}.
	\]
	we can therefore conclude that the aforementioned proof can be reread with $\op_W \derd{U^T}{u}$ and still applies, if only vectors $v \in W$ are allowed, or equivalently $\op_W \derd {U^T} u$ is entered instead of $\derd  {U^T} u$. This in turn equals $s$.
	\end{proof}
\end{lemma}
\begin{remark}
	Abgrall and collaborators in \cite{abgrall2019reinterpretation} showed that the entropy correction terms derived for the residual distribution schemes to enforce entropy conservation can be interpreted as solutions to optimization problems. These correction terms or the respective solutions to the optimization problems lead to the same descent directions, albeit their optimization problems are different. A closer inspection reveals that the optimization problems are in fact more or less dual to ours. Still, we are not interested in entropy conservation but high dissipation restricted by error bounds.
\end{remark}

Knowledge of the solution also allows us to give a lower bound for $\delta$ that makes our scheme even entropy dissipative in the classical sense.

\begin{theorem}[Classical entropy inequality for DDG]
	Assume that a monotone entropy stable FV flux $f(u_l, u_r)$ and a strictly convex and twice continuously differentiable entropy $U$ is used. Let $\epsilon$ be determined by the following expression, depending on $u$
	\begin{equation}
		\begin{aligned}
		\epsilon \geq \lim_{N\to \infty}\frac{\norm{h}_{T, 2}}{\skp{h}{\derd{U^T}{u}}_T}\left(\delta(N)\norm{\widetilde{\derd {U^T} u}}_{T, 2}  + \delta_U \norm{\derive{u^{T, N}}{t}}_{T, 1}\right), \\
		 \delta(N) = \norm{\derd{u^T}{t} - \derive{u^{\mathrm{ref, N}}}{t}}_{T, 2}, \quad
		   \delta_U = \norm{\derd{U}{u}\left(u^T\right) - \derd{U^T}{u}}_{T, \infty}.
		   \end{aligned}
		\label{eq:ess}
	\end{equation}
	Then hold the two entropy inequalities
	\[
	\derd{E_{u, T}} t \leq {F^T_l - F^T_r}, \quad \derd{E_{u, T}}{t}(t_0) \leq \lim_{N \to \infty}\derd{~}{t}\int_T U\left(u^{T, N}\right)\intd x,
	\]
	for our modified scheme
	\[
		\derd{u^{T, D}}{t} = \derd{u^T}{t} - s.
	\]
	The first inequality states that our solution satisfies a discrete entropy inequality, while the second one states that the entropy decreases faster than the entropy of the limit solution of the subcell scheme.
	\begin{proof}
		Let us remark before we start with our proof, that our corrections to the base scheme have zero mean value, i.e. the mean value of the time derivative of our scheme and of the base scheme are the same. We will exploit this behavior in what follows and use 
		\[
			\derd{u^T}{t} = \overline{\derd{u^T}{t}} + \widetilde{\derd{u^T}{t}}
		\]
		as a notation for the splitting of a function $g$ into a constant function $\overline g$ with the mean values of $g$ and a function $\tilde g$ representing the variation of $g$ around it's mean value.
		If $u^T = \overline {u^T}$ holds we can look at the special case of our base FV scheme as in this case $\derd{U}{u}(u(x_k)) = \derd{U}{u}(\bar u)$
			\[
				\derd{E_{u, T}(u(t))}{t} = \skp{\derd{U}{u}}{\derd{u^T}{t}}_T = \derd{U}{u}(\bar u) \skp{\One}{\derd{u^T}{t}}_T = \derd{U}{u}(\bar u)(f_l^T-f_r^T) \leq (F^T_l - F^T_r)
			\]
			is satisfied.
			Here we used that the entropy flux $F(u_l, u_r)$ of the entropy stable flux $f(u_l, u_r)$ satisfies
			\[
				\derd U u (u_m)(f(u_l, u_m) - f(u_m, u_r)) = \derd U u \derd {u_m} t = \derd U t \leq F(u_l, u_m) - F(u_m, u_r).
			\]
			See for example \cite{Tadmor87, Tad2003}.
			This entropy inequality also allows us to proof 
			\begin{equation}
			\lim_{N \to \infty} \skp{\derd U u (u^{T, N})}{\derd{u^{T, N}}{t}} =	\lim_{N \to \infty} \derd{~}{t}\int_T U(u^{T, N}(x, t))\intd x \leq {F^T_l - F^T_r}
				\label{eq:fvsei}
			\end{equation} for arbitrary $u^T(x, t)$,
			as this holds for the summed contributions of the FV subcell scheme in the cell for fixed $N$, and also in the limit. Clearly holds also
			\[
				\lim_{N \to \infty} \skp{\derd{U}{u}(u^T)}{ \derive{u^{T, N}}{t}} = \lim_{N \to \infty} \skp{\derd{U}{u}\left(u^{T, N}\right)}{\derive{u^{T, N}}{t}}
			\]
			as $\norm{\derive{u^{T, N}}{t}}_{T, 1}$ stays bounded while $\norm{\derd{U} u \left(u^T\right)- \derd U u \left(u^{T, N}\right)}_{T, \infty} \xrightarrow{N \to \infty} 0$ holds as the ansatz is continuous in the cell.  
		We will assume from now on $u^T \neq \bar u$ and therefore also that 
		\[
			h = g - \frac{\skp {\One}{g}}{\skp{\One}{\One}} \One \neq 0
		\]
		because $\derd{U}{u}(u^T) = 0$ is only possible for a single $u \in \R^n$ and $u \neq \bar u$ implies $\derive{U}{u} \neq \overline{\derive{U}{u}}$ as the entropy is strictly convex. We will now generalize our argument from the last section concerning the entropy dissipativity of our approximate time derivative $\derd{u^T}{t}$. It is imperative to first concentrate on the case where $N$ is finite as done before for the entropy inequality for the subcell scheme and derive appropriate bounds. We will afterwards go over to the limit to prove the theorem. Let us denote by $\derd U u (u^T(x, t))$ the exact value of the entropy variables associated with the numerical solution $u^T(t)$ while 
		\[
			\derd {U^T} u (u^T (x, t)) = \sum_{j \in T} \phi_j(x) \derd U u(u^T(x_j, t)), \text{where } \phi_j(x) \text{ are Lagrange polynomials},
		\]
		shall be the interpolation of $\derd U u (u^T(x, t))$ in the space $V$. We can therefore state the error made in the prediction of the entropy dissipation by interchanging the exact entropy functional with the one living in our approximation space
		\[
			\begin{aligned}
			&\skp{\derd U u \left(u^T(\cdot, t)\right)}{\derive{u^{T, N}}{t}}_T \\
			&=  \skp{\derd {U^T} u}{\derive{u^{T, N}}{t}}_T + \underbrace{\skp{\derd U u \left(u^T(\cdot,  t)\right) - \derd {U^T}{u}(\cdot, t)}{\derive {u^{T, N}} t}_T }_{\leq \delta_U \norm{\derive {u^{T, N}} t}_{T, 1}}.
			\end{aligned}
		\]
		Because $\derd{U^T}{u}$ is from $V$ we can exchange $\derive{u^{T, N}}{t}$ for $\derive{u^{\mathrm{ref}, N}}{t}$
		\[
			\skp{\derd {U^T}{u}}{\derive{u^{T, N}}{t}} = \skp{\derd {U^T} u }{\derive{u^{\mathrm{ref, N}}}{t}}
		\]
		without any penalty as $\derive{u^{T, N}}{t} - \derive{u^{\mathrm{ref, N}}}{t} \in V^\perp$ holds. 
		If we also swap the reference time derivative for our scheme derivative $\derd{u^{T, N}}{t}$, we find
		\[
		\begin{aligned}
			&\skp{\derd{U^T}{t}}{\derive{u^{\mathrm{ref, N}}}{t}}_T 
			\\ =& \skp{\overline{\derd{U^T}{u}} + \widetilde{\derd{U^T}{u}}}{ \overline{\derive{u^{\mathrm{ref, N}}}{t}} + \widetilde{\derive{u^{\mathrm{ref, N}}}{t}}}_T
			\\ =&
			\underbrace{\skp{\overline{\derd{U^T}{u}} + \widetilde{\derd{U^T}{u}}}{\overline{\derive{u^{\mathrm{ref, N}}}{t}}}_T}_{\text{same mean value}} + 
			\underbrace{\skp{\overline{\derd{U^T}{u}}}{\widetilde{\derive{u^{\mathrm{ref, N}}}{t}}}_T}_{\text{orthogonality}}
			+
			\skp{\widetilde{\derd{U^T}{u}}}{\widetilde{\derive{u^{\mathrm{ref, N}}}{t}}}_T
			\\ =&
			\underbrace{\skp{\overline{\derd{U^T}{u}} + \widetilde{\derd{U^T}{u}}}{\overline{\derd{u^{T}}{t}}}_T +
			\skp{\overline{\derd{U^T}{u}}}{\widetilde{\derd{u^{T}}{t}}}_T + 
			\skp{\widetilde{\derd{U^T}{u}}}{\widetilde{\derd{u^{T}}{t}}}_T}_{\skp{\derd{U^T}{u}}{\derd{u^{T}}{t}}_T}\\
			&+
			 \underbrace{\skp{\widetilde{\derd{U^T}{u}}}{\widetilde{\derive{u^{\mathrm{ref}, N}}{t}} - \widetilde{ \derd{u^{T}}{t}}}_T}_{\leq \norm{\widetilde{\derd {U^T} u}}_{T, 2}\delta(N) }.
			\end{aligned}
		\]
		Please note that we used several facts of our reference solution and our approximate solution to sharpen this bound to only depend on the variation of the entropy variables and the variation of the solution around their respective mean values. This was possible after the inner product was split into the respective inner products of the mean values and variations around the mean values with each other. As our reference solution $\derive{u^{\mathrm{ref}}}{t}$ and $\derd{u^{T}}{t}$ share the same mean values we can easily swap one for the other in the first inner product. The second inner product in this decomposition can also be swapped, as it is by definition zero. This follows from the fact that the mean values, calculated with respect to the inner product of $T$, are orthogonal to the variations. The only penalty that has to be bounded consists of the inner product of the variations.
		Combining the previous steps leads us to
		\[
		\begin{aligned}
			&\abs{\skp{\derd U u \left(u^T(t)\right)}{\derive {u^{T, N}} t}_T - \skp{\derd {U^T} u (t)}{\derd {u^T} t}_T} &&\\
			\leq& 
			 \abs{\skp{\widetilde{\derd{U^T} u (t)}}{\widetilde{\derive{u^{\mathrm{ref}, N}}{t} }-\widetilde{\derd{u^T}{t}}}_T} &+& \abs{\skp{\derd U u - \derd {U^T} u}{ \derive {u^{T, N}} t}_T} \\
			 \leq&
			 \norm{\widetilde{\derd {U^T} u}}_{T, 2} \delta(N)  &+& \delta_U \norm{\derive {u^{T, N}} t}_{T, 1}
			 \end{aligned}
		\]
		as an upper bound for the difference of the exact entropy dissipation and the entropy dissipation of our approximate solution. If $\epsilon$ is set according to the value given above, one finds that the following inequality holds as the additional entropy production that can be bounded using the error $\delta$ and $\delta_U$ can be fully counteracted by the entropy dissipation of the steepest descent direction
		\[
		\begin{aligned}
		\skp{\derd{U^T}{u}}{\derd {u^T} t - s_N}_T 
		=
		\skp{\derd{U^T}{u}}{\derd {u^T} t - \epsilon_N \frac{h}{\norm h}}_T 
		=
		\skp{\derd {U^T} u}{\derd {u^T} t}_T - \epsilon \frac{\skp{h}{\derd{U^T}{u}}_T}{\norm{h}_{T, 2}} \\
		 = 
		 \skp{\derd {U^T} u}{\derd {u^T} t}_T - \left(\delta \norm{\widetilde{\derd {U^T} u}}_{T, 2} + \delta_U\norm{\derive {u^{T, N}} t}_{T, 1}\right) \leq \skp{ {\derd U u\left (u^T(t)\right)}}{ {\derive {u^{T, N}} t}}.
		\end{aligned}
		\]
		The given lower bound for $\epsilon$  behaves significantly better than one would think. A primary reason for this is that $h$ is co-linear to $\widetilde{ \derd{U^T}{u}}$ and relates to $\derd{U^T}{u}$ via an orthogonal projection. Therefore is the expression 
		\[
			\frac{\norm h _{T, 2} \norm{\widetilde{ \derd{U^T}{u}}}_{T, 2}}{\skp{h}{\derd{U^T}{u}}_T} = \frac{\norm h _{T, 2} \norm{\widetilde{ \derd{U^T}{u}}}_{T, 2}}{\skp{h}{\widetilde{\derd{U^T}{u}}}_T +\skp{h}{\overline{\derd{U^T}{u}}}_T }
		\]
		bounded.
		We showed earlier, that the limit of the reference solution exists. We can therefore go over to the limit and conclude that
		\[
			\skp{\derd{U^T}{u}}{\derd{u^{T, D}}{t}}_T = \skp{\derd{U^T}{u}}{\derd {u^T} t - s}_T \leq  \lim_{N \to \infty}\derd{~}{t}\int_T U\left(u^{T, N}\right)\intd x
		\]
		holds in the limit.
		The last step in our proof is a comparison principle. We combine the last equation with \eqref{eq:fvsei} to find
		\[
			\skp{\derd{U^T}{u}}{\derd{u^{T, D}}{t}}_T = \skp{\derd{U^T}{u}}{\derd u t - s}_T \leq F^T_l - F^T_r.
		\]
	\end{proof}

\end{theorem}

	\section{Construction of the discrete scheme} \label{sec:disccons}
	Designing an error estimator for the discrete scheme is a more delicate issue. The observation
\[
\begin{aligned}
\derd{~}{t}\norm{u^{\refs}(\cdot, t) - u^T(\cdot, t)}_{T, 2}
& = 
 \derd{~}{t}\sqrt{\skp{u^{\refs}(\cdot, t) - u^T(\cdot, t)}{u^{\refs}(\cdot, t) - u^T(\cdot, t)}_{T}}\\
&= \frac{\derd{~}{t} \skp{u^{\refs}(\cdot, t) - u^T(\cdot, t)}{u^{\refs}(\cdot, t) - u^T(\cdot, t)}_{T} }{2 \sqrt{\skp{u^{\refs}(\cdot, t) - u^T(\cdot, t)}{u^{\refs}(\cdot, t) - u^T(\cdot, t)}_{T}}} \\
&= 
 \frac{\skp{u^\refs - u^T}{\derive{u^\refs}{t} - \derd{u^T}{t}}_T}{\norm{u^\refs - u^T}_{T, 2}} 
  \leq \norm{\derive{u^\refs}{t} - \derd{u^T}{t}}_{T, 2} = \delta^T_t
 \end{aligned}
\]
allows us to estimate the total error made in cell $T$ if we know the error between the derivative of the reference solution and the exact solution,
\[
	\norm{u^{\mathrm{ref}}(x, t) - u^T(x, t)}_{T, 2} \leq \int_0^t \norm{\derive{u^\refs}{t}(\tau) - \derd{u^T}{t}(\tau)}_{T, 2} \intd \tau = \int_0^t \delta_\tau^T \intd \tau.
\] But during the integration of such an error estimate in time the two solutions will in general start to drift apart from each other. Calculating a reference derivative would therefore need the knowledge of the exact solution. We will assume therefore that the total error between reference solution and numerical solution stays small enough to warrant us using the reference derivative at time $\tau$ as calculated from the solution $u^T(\cdot, \tau)$ and not with respect to $u^{\mathrm{ref}}(\cdot, \tau)$. 
As in the other cases, we will reside to numerical quadrature for the calculation of this quantity. The integrand was already used as an error estimate in the semidiscrete case. For the calculation of the outer integral in the time direction it is worthwhile to consider the connection between Runge-Kutta time integration and numerical quadrature. If for example the SSPRK33 solver \cite{SO1988}
\[
	\begin{aligned}
	u^{(1)} &= u^{(0)} + \Delta t L(u^{(0)}) \\
	u^{(2)} &= u^{(0)} + \frac 1 4 \Delta t L(u^{(0)}) + \frac 1 4 \Delta t L(u^{(1)})\\
	u^{(3)} &= u^{(0)} + \frac 1 6 \Delta t L(u^{(0)}) + \frac 1 6 \Delta t L(u^{(1)} )+ \frac 2 3 \Delta t L(u^{(2)})
	\end{aligned}
\] is used one sees clearly, that this in fact a numerical quadrature of $\derd{u}{t} = L(u)$ where first using the left sided Newton-Cotes formula an approximation $u^{(1)}$ for the solution at $\Delta t$ is calculated. Next an approximation $u^{(2)}$ for the solution at time $\Delta t /2$ is calculated using the trapezoidal rule
\[
	\int_0^{\Delta t} \derd{u}{t} \intd t \approx \frac{L(u(0)) + 4 L(u(\Delta t/2)) + L(u(\Delta t))}{6},
\] and as a last step a better approximation of the rightmost value is calculated using the Simpson rule and the two precalculated approximations. If our semi-discrete error estimate is local Lipschitz continouous, which it clearly is, it is therefore a logical decision to choose the quadrature of the time integrator as quadrature for discrete error estimate. One can therefore also reuse the interim results of the time integrator. Our discrete error estimate therefore reads as
\begin{equation}
	\delta^T = \Delta t \frac{\delta^T_{t}\left(u^{(0)}\right) + 4 \delta^T_t\left(u^{(2)}\right) + \delta^T_t\left(u^{(1)}\right)}{6}.
	\label{eq:tserror}
\end{equation}
and is third order accurate, as is the timesteping algorithm.
\begin{lemma}[Existence of at most one solution]
	The optimization problem stated for the DRKDG method \eqref{eq:discprob} possesses a unique solution if a strictly convex entropy functional is used.
	\begin{proof}
		The set $Z$ is closed and bounded, hence compact, and the functional $E^T(u)$ continuous. Therefore there exists a $u$ in $Z$ with $\forall v \in Z: E^T(u) \leq E^T(v)$.
		Concerning the uniqueness, the strict convexity of $E^T(\cdot)$ implies that if $u, v \in V$, with $u \neq v$, would exist with $E^T(u) = E^T(v)$ and $\forall w \in W: E^T(u) \leq E^T(w)$ a contradiction for $\lambda = \frac 1 2$ and $w = \lambda u + (1-\lambda) v$ would arise as
		\[
		E^T(w) = E^T(\lambda u + (1-\lambda) v) < \lambda E^T(u) + (1-\lambda) E^T(v) = E^T(u)
		\] 
		would follow from the strict convexity in this case.
	\end{proof}
\end{lemma}
Solving the optimization problem for the discrete case is more intricate than in the semidiscrete case. While a solution for the simple entropy $U(u) = \frac{u^2}{2}$ could be computed by hand, the computation for complicated entropies is not feasible. We will therefore concentrate on the numerical approximation of a solution. A simple, yet effective, procedure seems to be the gradient descent
\[
u^{T, n+1} = u^{T, n} - \lambda h
\]
with appropriately chosen step size $\lambda$ \cite{Nesterov2004Convex}. The descent direction will be, as in the semidiscrete case, the solution to the steepest descent problem from Lemma \ref{lem:rgd}. As we would like to solve this problem in every complete RK step should the cost of the nonlinear solver be of the same magnitude as the repeated evaluations of the semidiscretisation for the Runge-Kutta time integration method. This limits the number of allowed steps and we generally try to use the same amount of gradient steps $r$ as there are stages in the RK method. We therefore propose 
\[
\lambda = \frac{\epsilon}{r \norm{h}_{T, 2}}
\]
as this limits the maximal distance between $u^{T, 0}$ and the final result $u^{T, s}$ to
\[
\norm{u^{T, 0} - u^{T, r}}_{T, 2} \leq \sum_{n = 1}^r \norm{u^{T, n} - u^{T, n-1}}_{T, 2} \leq \sum_{n=1}^r \frac{\epsilon \norm h_{T, 2}}{r \norm h_{T, 2}} = \epsilon.
\]
One restriction has to be made, as this selection could lead to a diverging sequence of steps if $\lambda$ is large as $\epsilon$ is large. We therefore restrict the step size to counteract this problem. Let $L$ be Lipschitz bound on the entropy variables in the norm on $T$. Then follows from \cite[equation 1.2.5]{Nesterov2004Convex}
\[
	E^T\left(u^{T, n+1}\right)\leq E^T\left(u^{T, n}\right) + \skp{\derd {U^T} u }{\lambda h}_T + \frac L 2 \lambda^2 =  E^T\left(u^{T, n}\right) - \lambda \underbrace{\skp{\derd {U^T} u }{-h}_T}_{\geq 0}  + \frac L 2 \lambda^2.
\]
This implies that a descent happens whenever the step size is chosen as
\[
\lambda < 2\frac{\skp{\derd{U^T}{u}}{-h}_T}{L}.
\]
In fact, in the implementation
\[
\lambda \leq \frac{3\skp{\derd{U^T}{u}}{-h}_T}{2L}
\]
is used as a compromise as
\[
\lambda \leq \frac{\skp{\derd{U^T}{u}}{-h}_T}{L}
\]
minimizes this bound.

	\section{Numerical tests} \label{sec:NT}
	Our numerical tests were carried out to shed some light on the following topics and possible problems:
\begin{itemize}
	\item Does the added error estimate controlled dissipation stabilize the schemes enough to use them for shock capturing calculations?
	\item Is the numerical Dafermos entropy rate criterion satisfied, even though only approximate errors for the norms are used?
	\item Are the added corrections small enough when a smooth solution is calculated to not destroy the high order approximation?
	\item Does the added dissipation reduce or enlarge the possible time step?
\end{itemize}
Before detailing the results, let us state that all of these questions can be answered with promising results, although the question of small approximation errors needs further research. The tested methods, i.e. their used quadrature rules, collocation points, time integration and similar design decisions are given in the table \ref{tab:sd}.
\begin{table}
	\begin{tabular}{c|c|c|c}
	Option	& DDG &  DRK-DG & Godunov \\
	\hline
	Time Integrator & SSPRK33 & SSPRK33 & SSPRK33\\
	Timestep & $0.5 = \frac{\Delta t}{(p^2 + 1)\Delta x}c_{max}$ & $0.5 = \frac{\Delta t}{(p^2 + 1)\Delta x}c_{max}$& $\frac{\Delta t}{\Delta x}c_{max} = 0.5$ \\
	Stepsize & from \eqref{eq:ess} & from \eqref{eq:tserror}& none \\
	Collocation points &$p+1$ pt.~Gauß-Lob.  & $p+1$ pt.~Gauß-Lob. & none \\
	Error estimate & $p+2$ pt.~Gauß-Leg. & $p+2$ pt.~Gauß-Leg. & none \\
	Polynomial order $p$ &  6 & 6 & piecewise const. \\
	\end{tabular}
	\caption{Overview of the schemes that were tested against each other in the numerical tests section. Their respective parameters concerning time integration are given together with the used quadratures and collocation points. The Godunov method is used to calculate reference solutions on a fine mesh.}
	\label{tab:sd}
\end{table}

\subsection{Calculation of discontinuous solutions}
The DRKDG and DDG schemes were first tested for Burgers' equation
\[
	\Omega = [0, 2), \quad f(u) = \frac{u^2}{2}, \quad U(u) = u^2
\]
on a periodic domain with the square entropy and the local Lax-Friedrichs flux. The tested initial conditions were
\[
	u_1(x, 0) = \sin(\pi x) + \frac 1 2, \quad u_2(x, 0) = \begin{cases}-x & x \in [0, 1) \\2-x & x \in [1, 2)  \end{cases}
\]
and it is well known that $u_1$ results in a discontinuous solution in finite time and crashes a RK-DG method without further stabilization quite easily. The second initial condition results in a rarefaction of the initial discontinuity. Because the sonic point of the flux is also part of this rarefaction we will be able to analyze the behavior of the scheme in this sometimes troublesome situation \cite{Tang2005Sonic}.
\begin{figure}
	\begin{subfigure}{0.49\textwidth}
		\includegraphics[width=\textwidth]{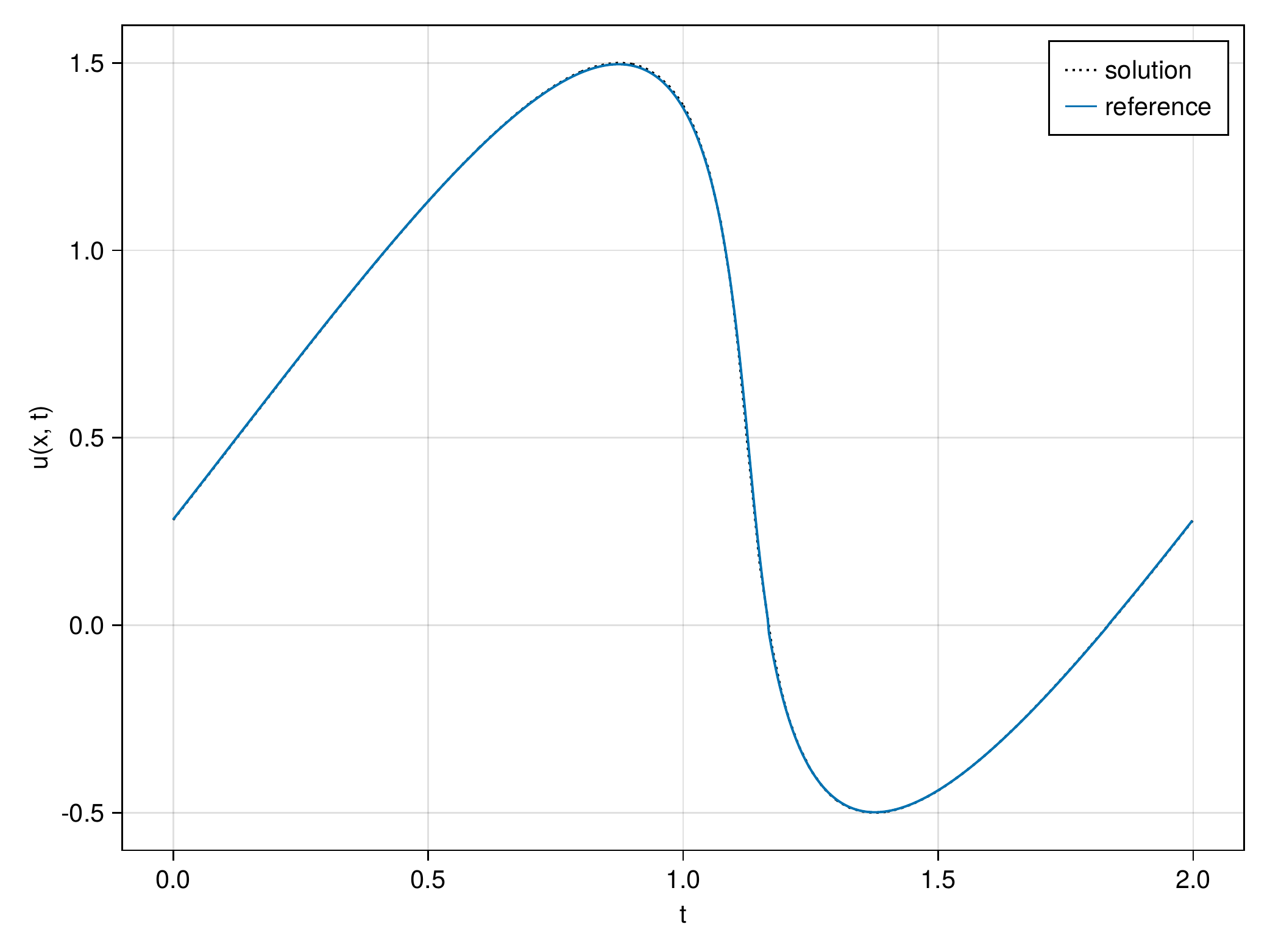}
		\caption{t = 0.25, 20 Elements}
	\end{subfigure}
	\begin{subfigure}{0.49\textwidth}
		\includegraphics[width=\textwidth]{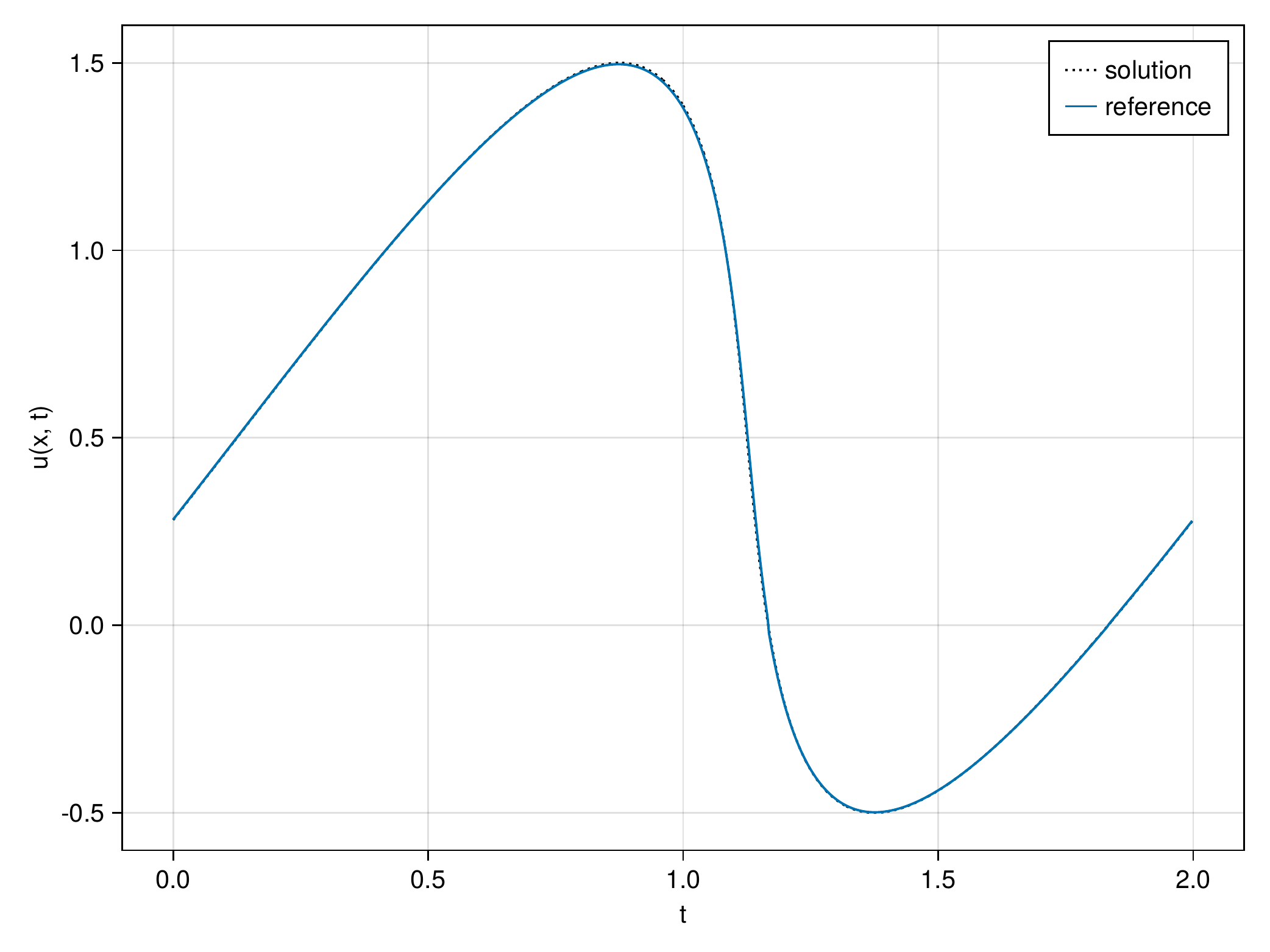}
		\caption{t = 0.25, 50 Elements}
	\end{subfigure}
	\begin{subfigure}{0.49\textwidth}
		\includegraphics[width=\textwidth]{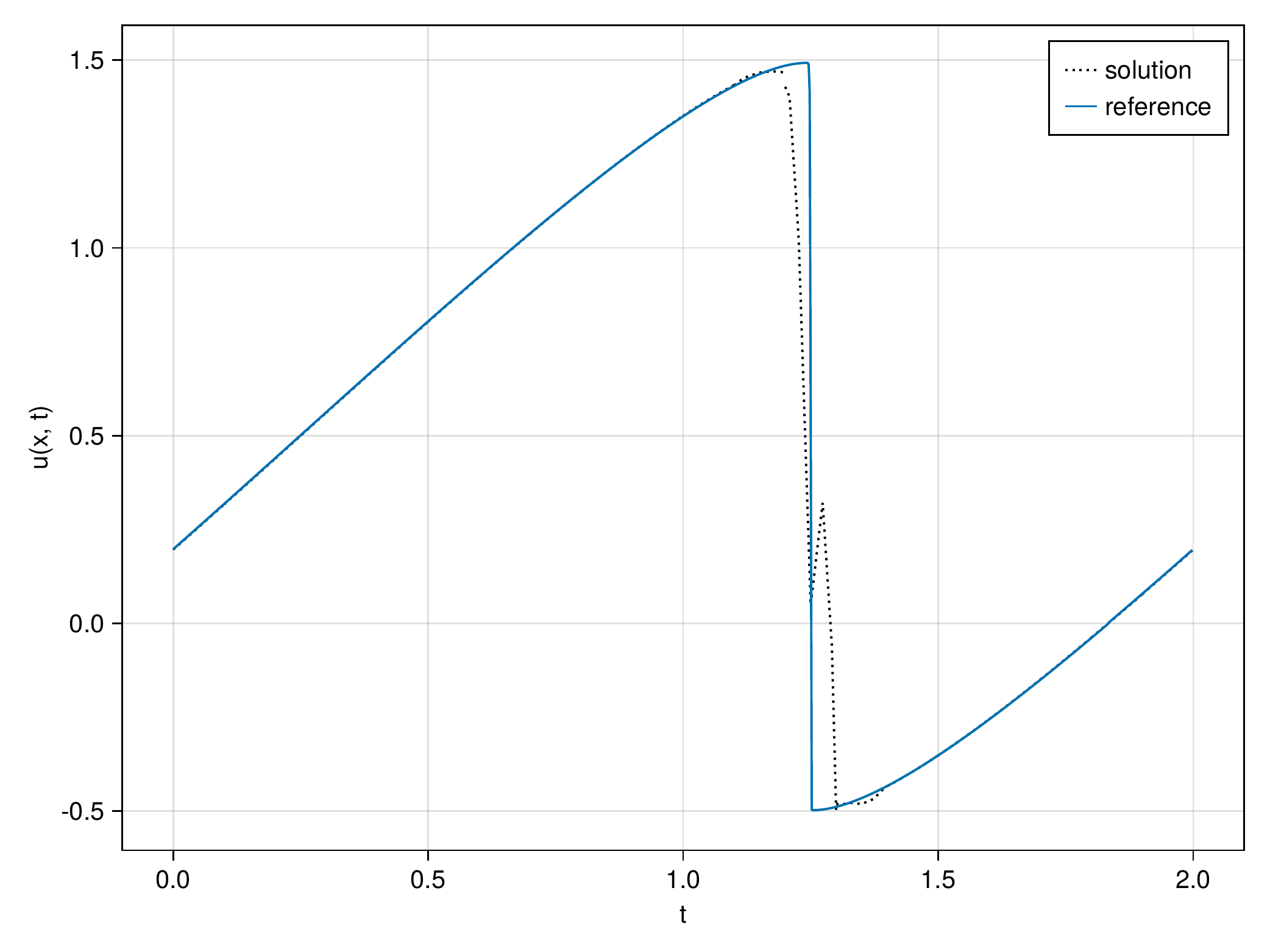}
		\caption{t = 0.5, 20 Elements}
	\end{subfigure}
	\begin{subfigure}{0.49\textwidth}
		\includegraphics[width=\textwidth]{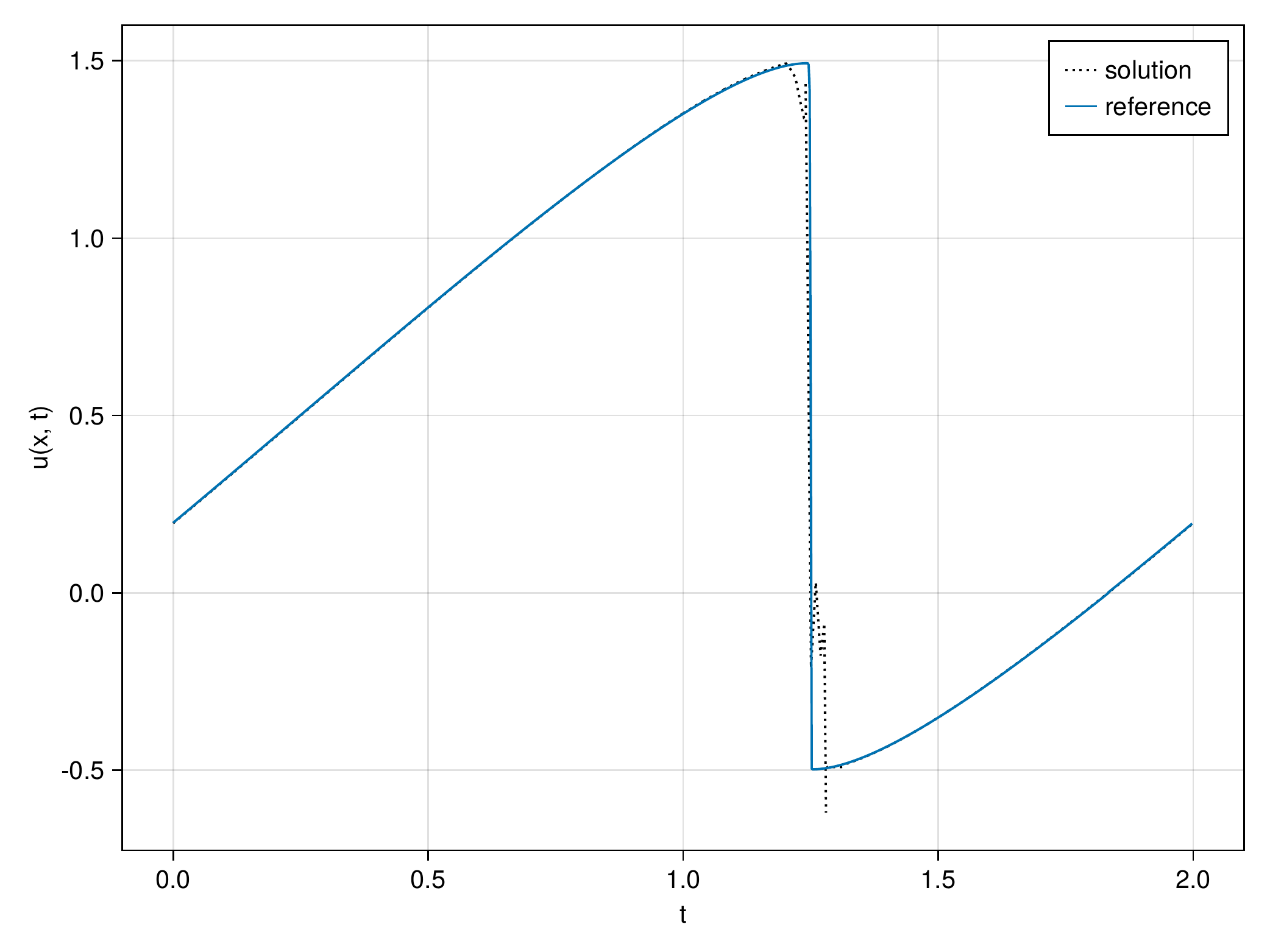}
		\caption{t = 0.5, 50 Elements}
	\end{subfigure}
	\begin{subfigure}{0.49\textwidth}
		\includegraphics[width=\textwidth]{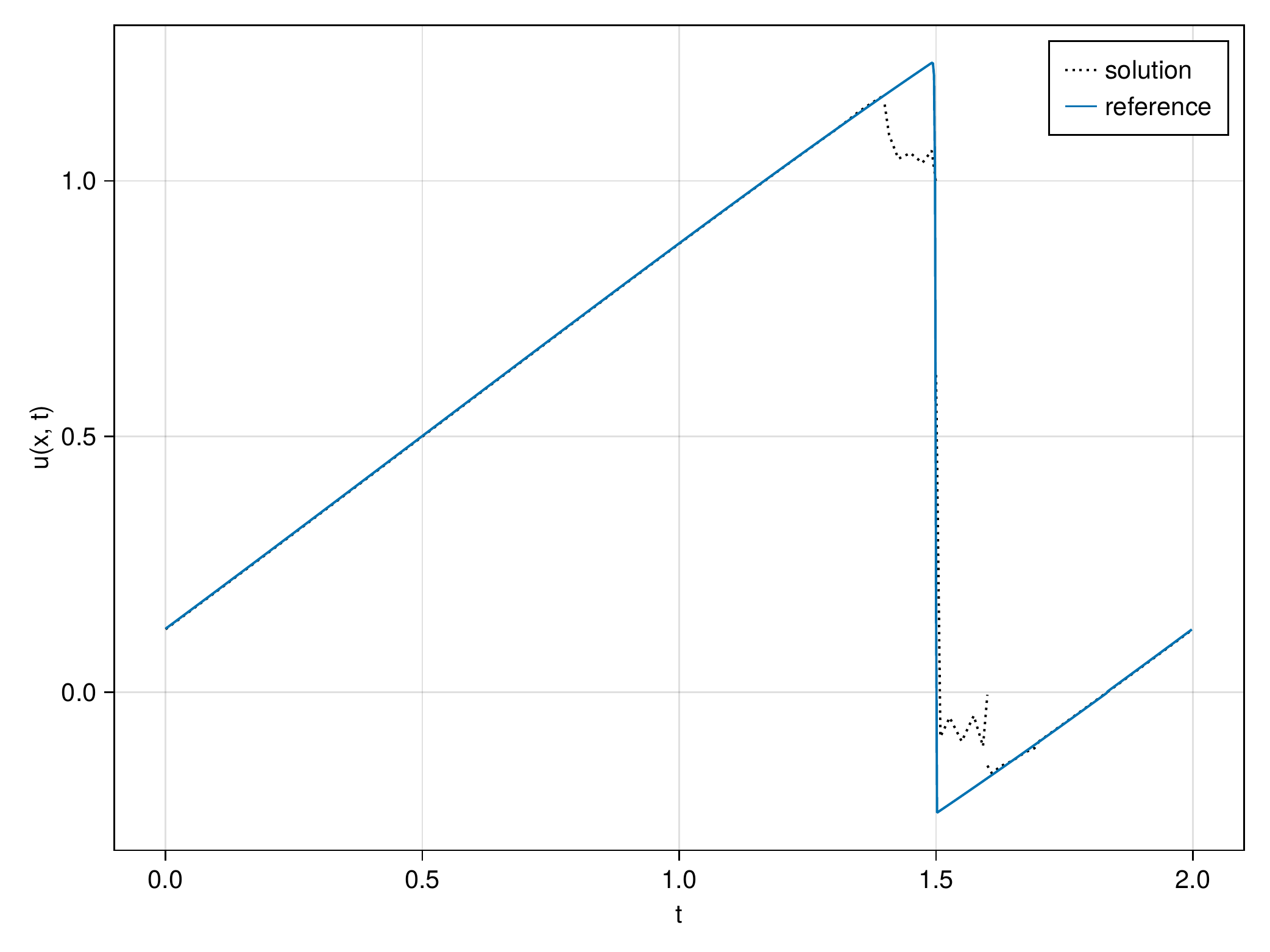}
		\caption{t = 1.0, 20 Elements}
	\end{subfigure}
	\begin{subfigure}{0.49\textwidth}
		\includegraphics[width=\textwidth]{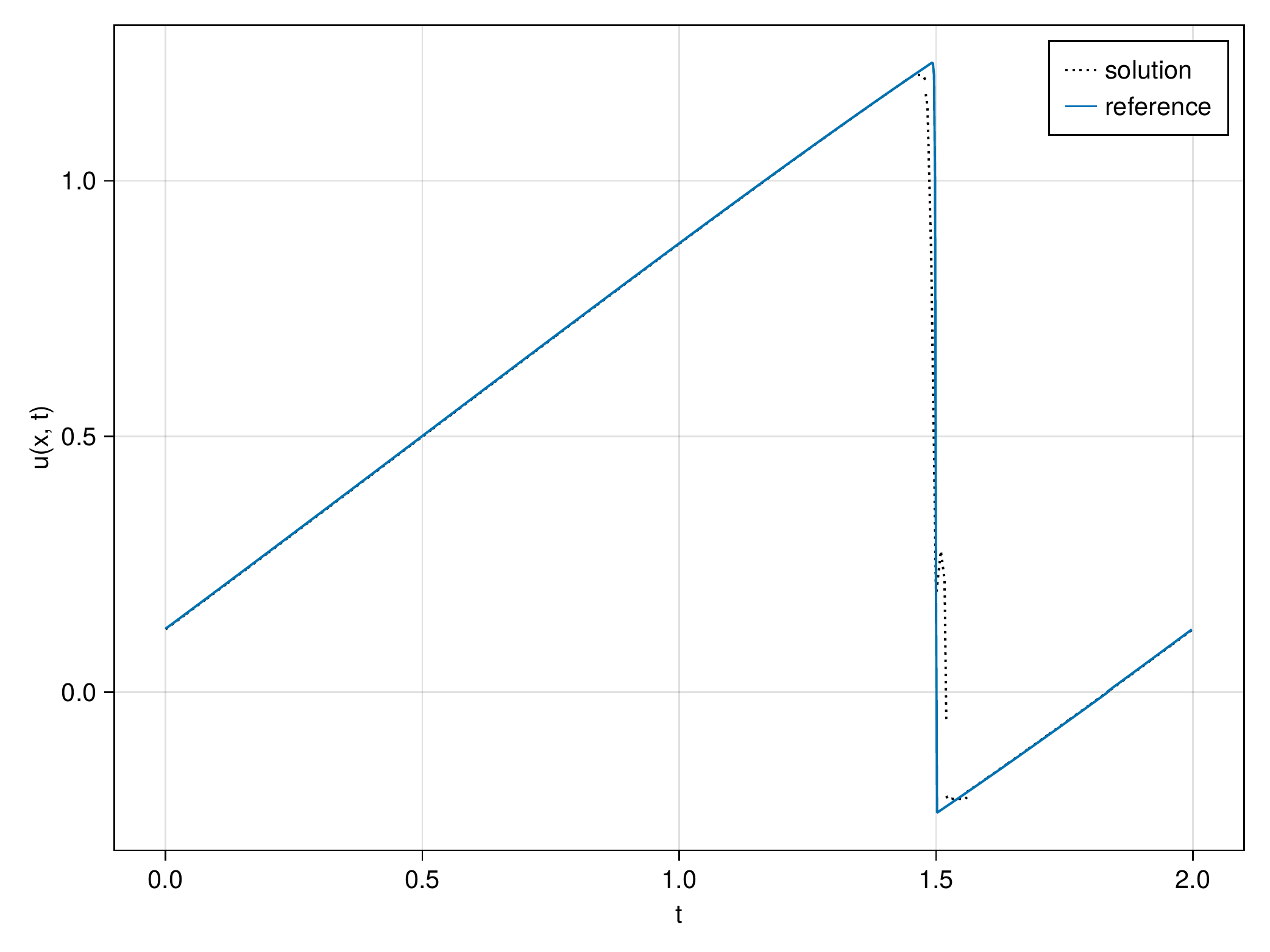}
		\caption{t = 1.0, 50 Elements}
	\end{subfigure}
	\caption{Solution to the first initial condition using the semi-discrete scheme DDG.}
	\label{fig:sintest}
\end{figure}
\begin{figure}
	\begin{subfigure}{0.49\textwidth}
		\includegraphics[width=\textwidth]{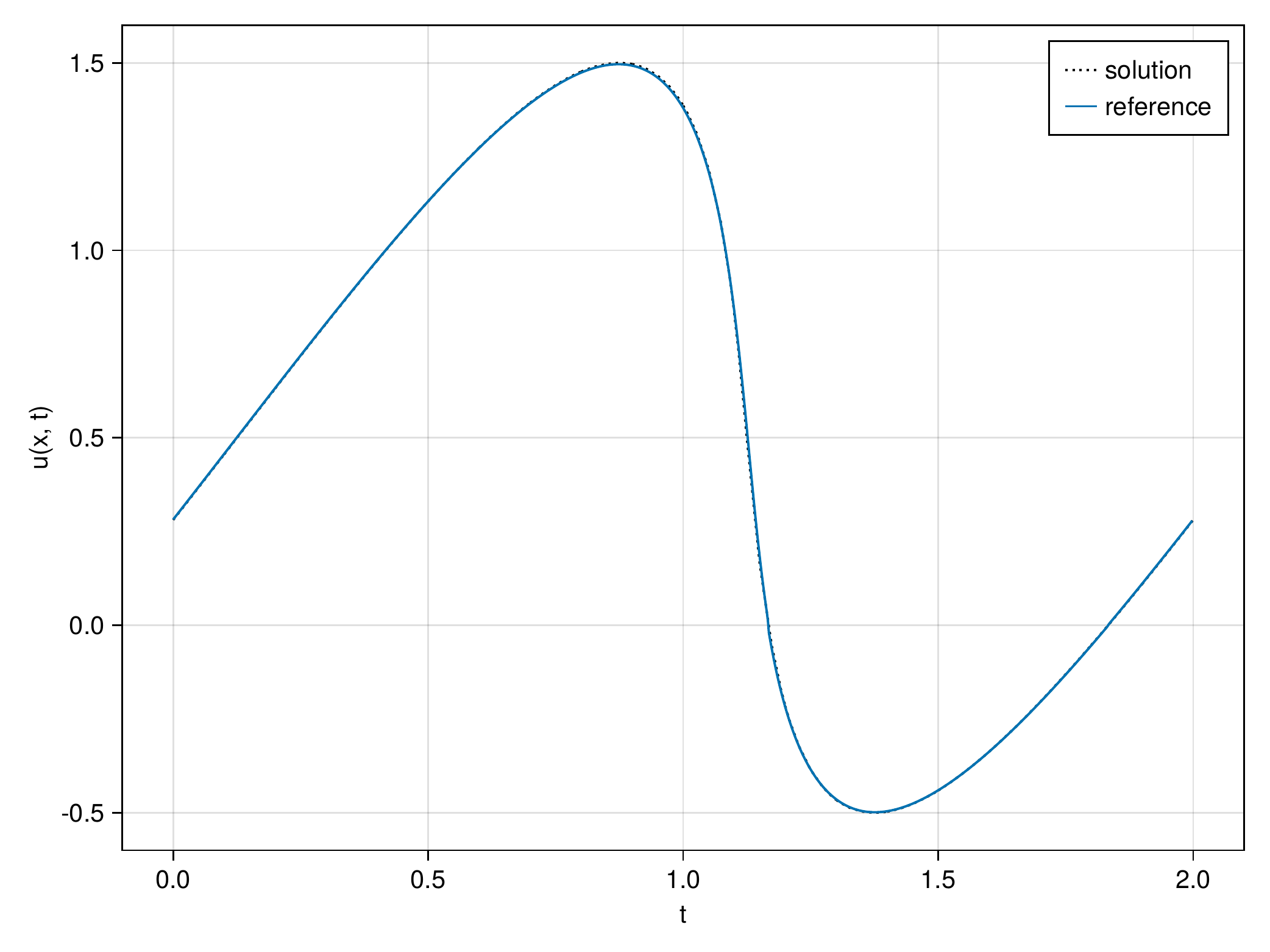}
		\caption{t = 0.25, 20 Elements}
	\end{subfigure}
	\begin{subfigure}{0.49\textwidth}
		\includegraphics[width=\textwidth]{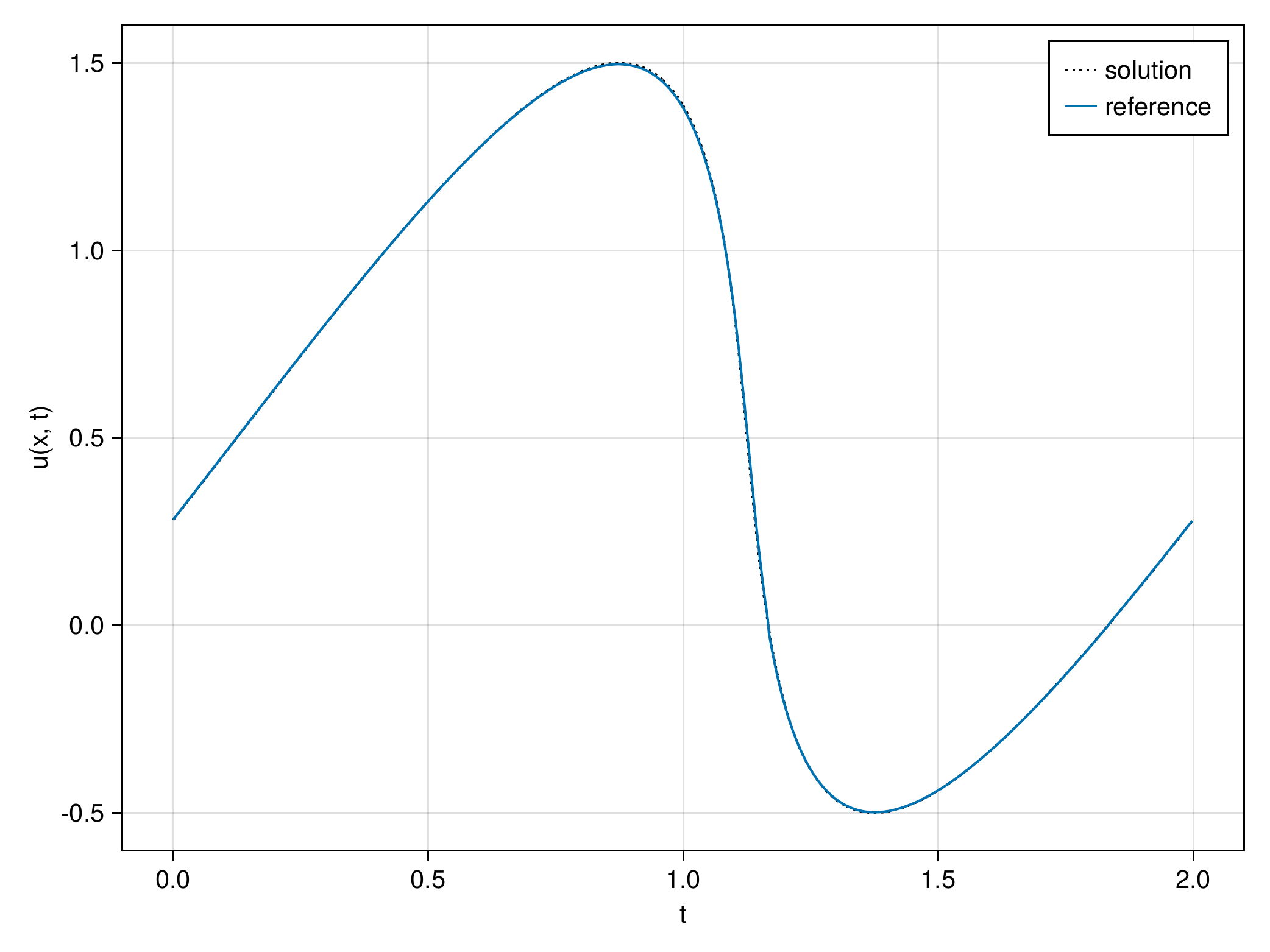}
		\caption{t = 0.25, 50 Elements}
	\end{subfigure}
	\begin{subfigure}{0.49\textwidth}
		\includegraphics[width=\textwidth]{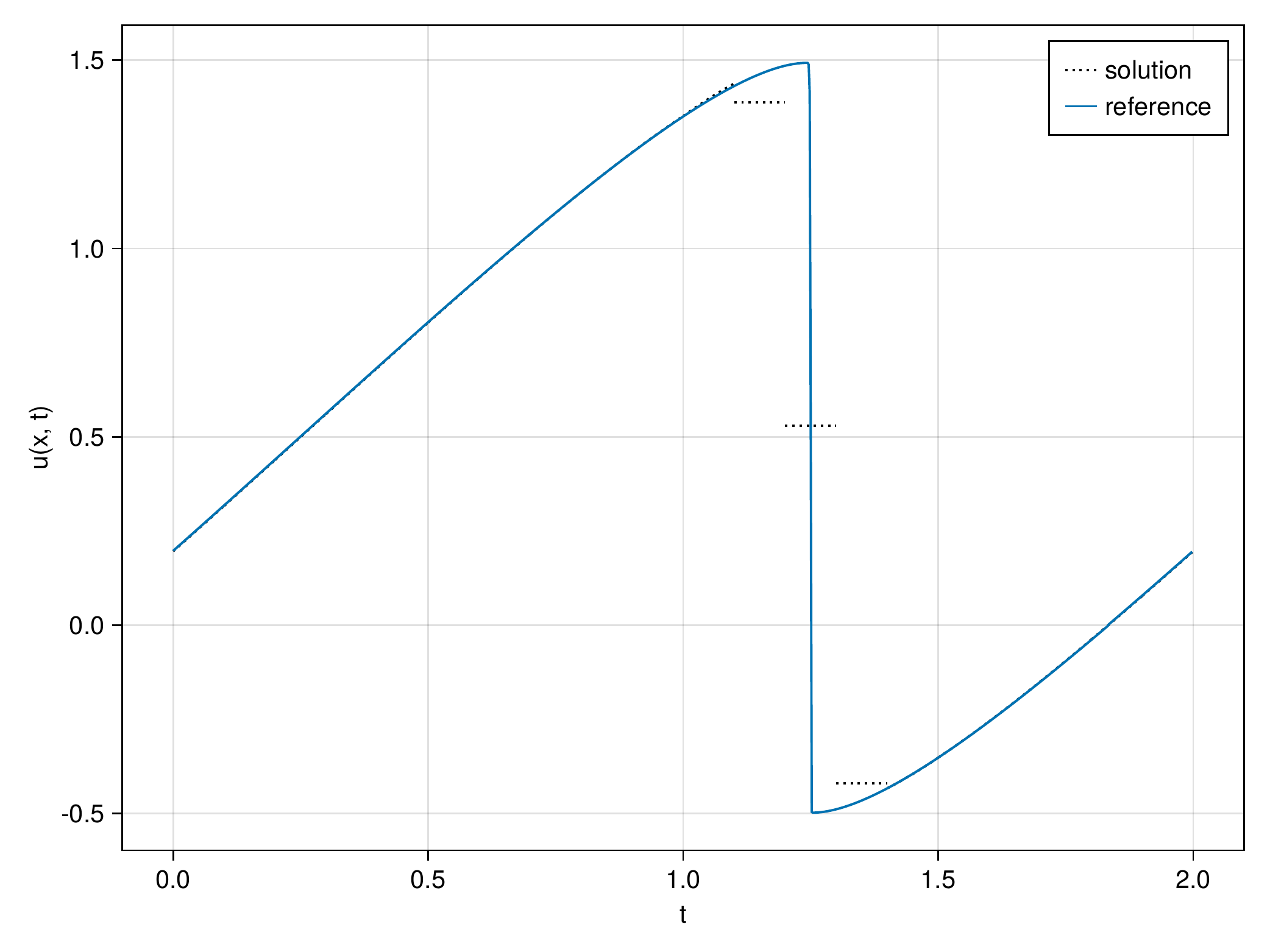}
		\caption{t = 0.5, 20 Elements}
	\end{subfigure}
	\begin{subfigure}{0.49\textwidth}
		\includegraphics[width=\textwidth]{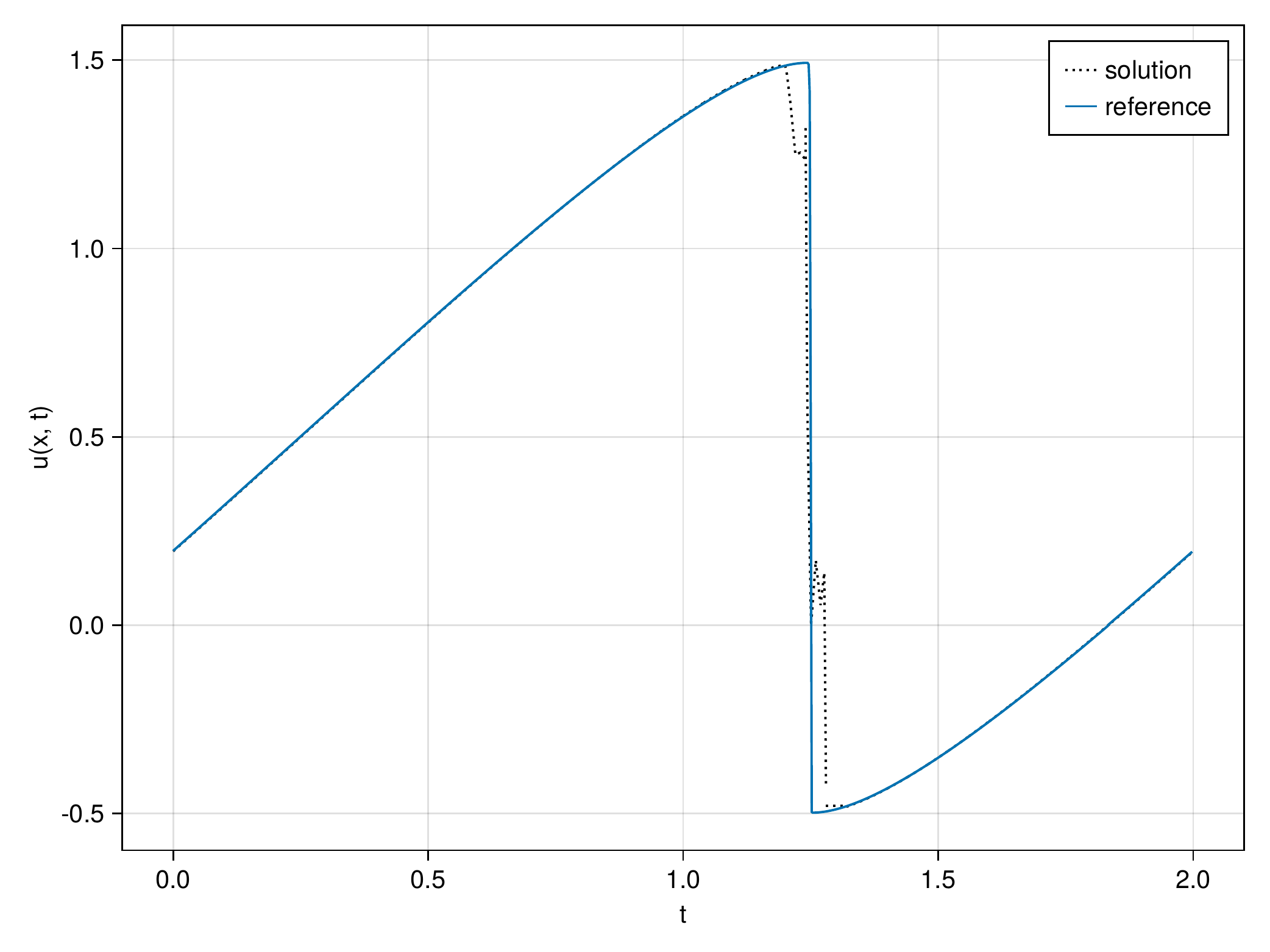}
		\caption{t = 0.5, 50 Elements}
	\end{subfigure}
	\begin{subfigure}{0.49\textwidth}
		\includegraphics[width=\textwidth]{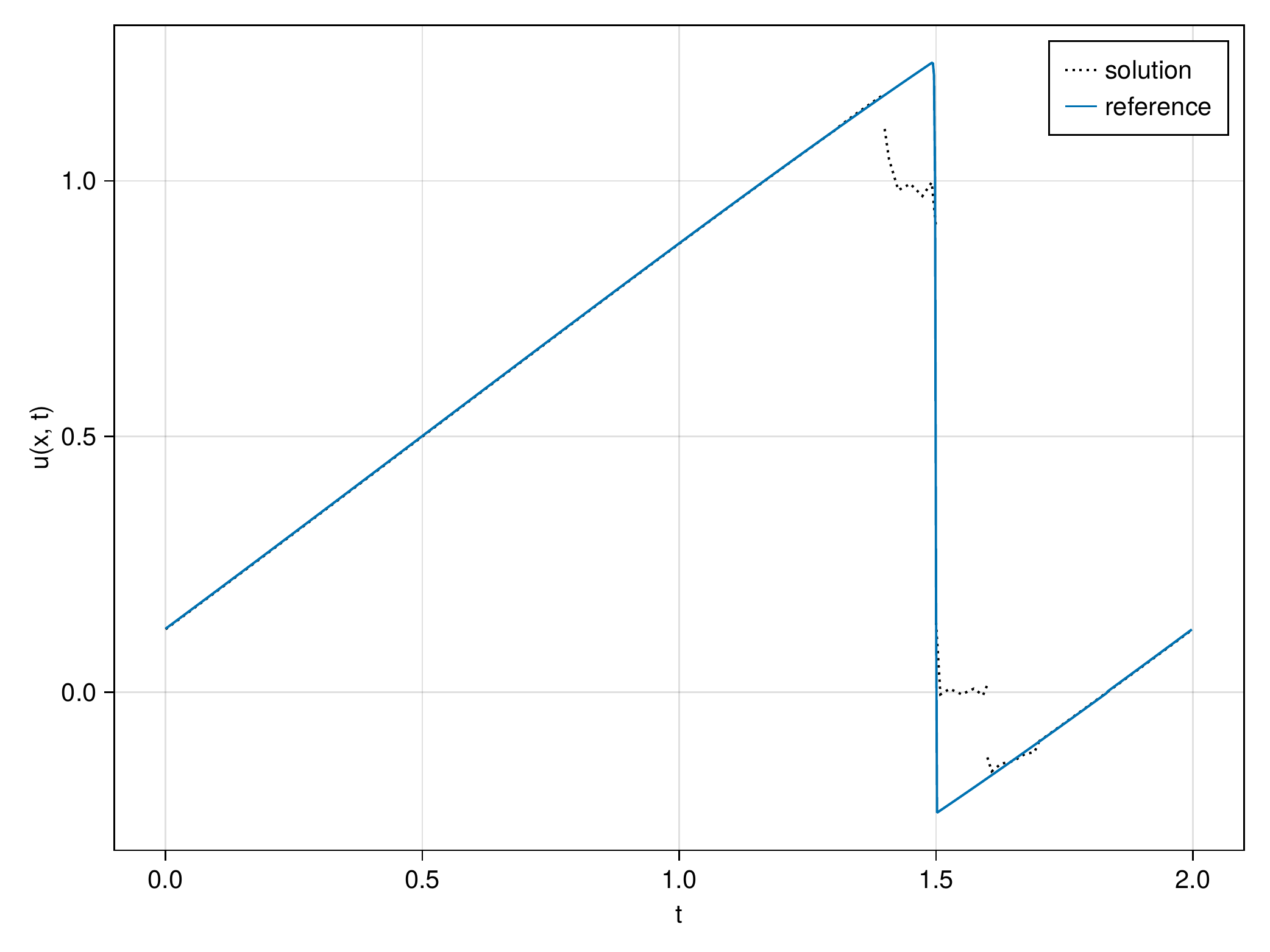}
		\caption{t = 1.0, 20 Elements}
	\end{subfigure}
	\begin{subfigure}{0.49\textwidth}
		\includegraphics[width=\textwidth]{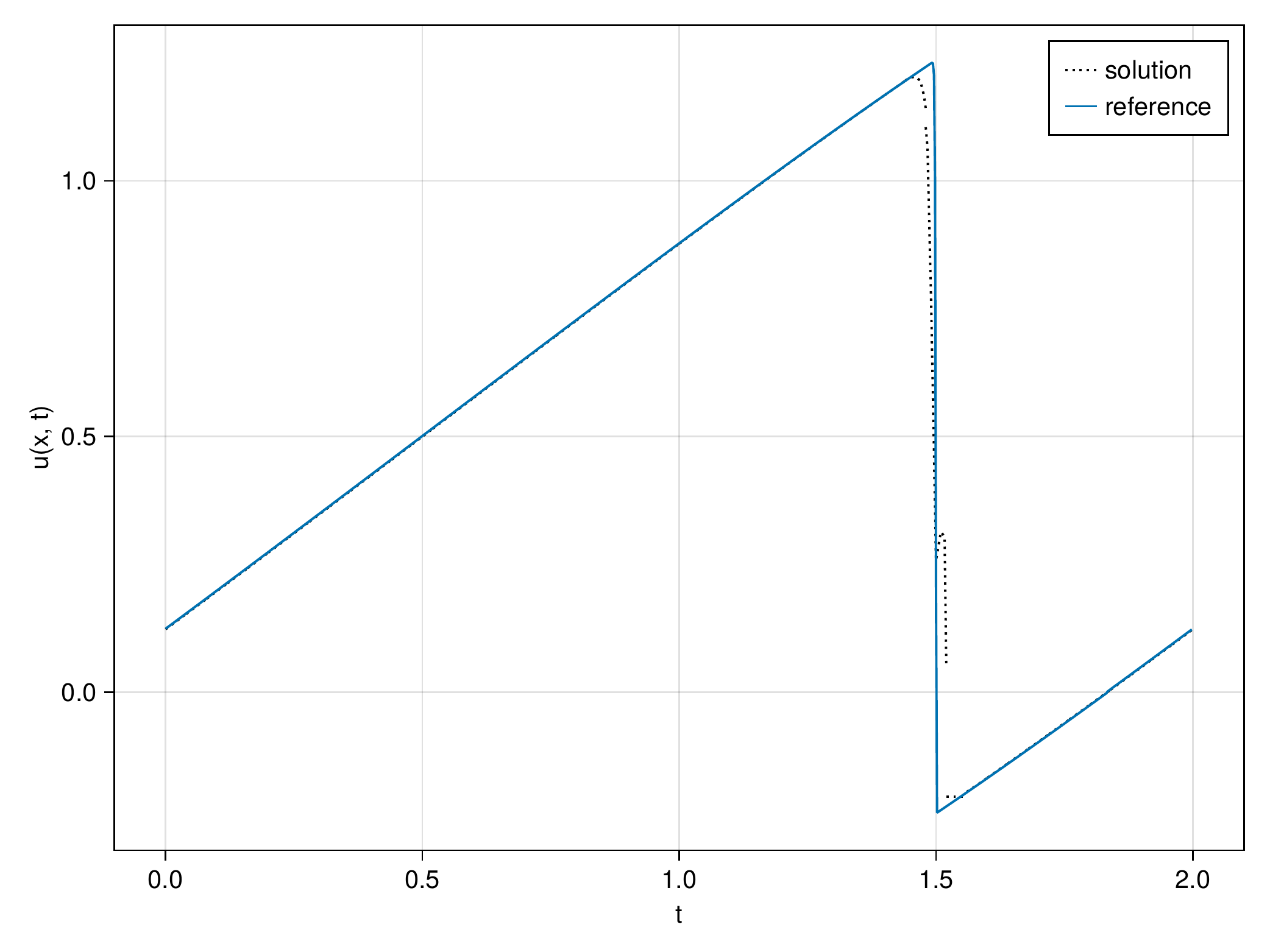}
		\caption{t = 1.0, 50 Elements}
	\end{subfigure}
	\caption{Solution to the first initial condition by the fully discrete scheme DRKDG.}
	\label{fig:dst}
\end{figure}
\begin{figure}
	\begin{subfigure}{0.49 \textwidth}
		\includegraphics[width=\textwidth]{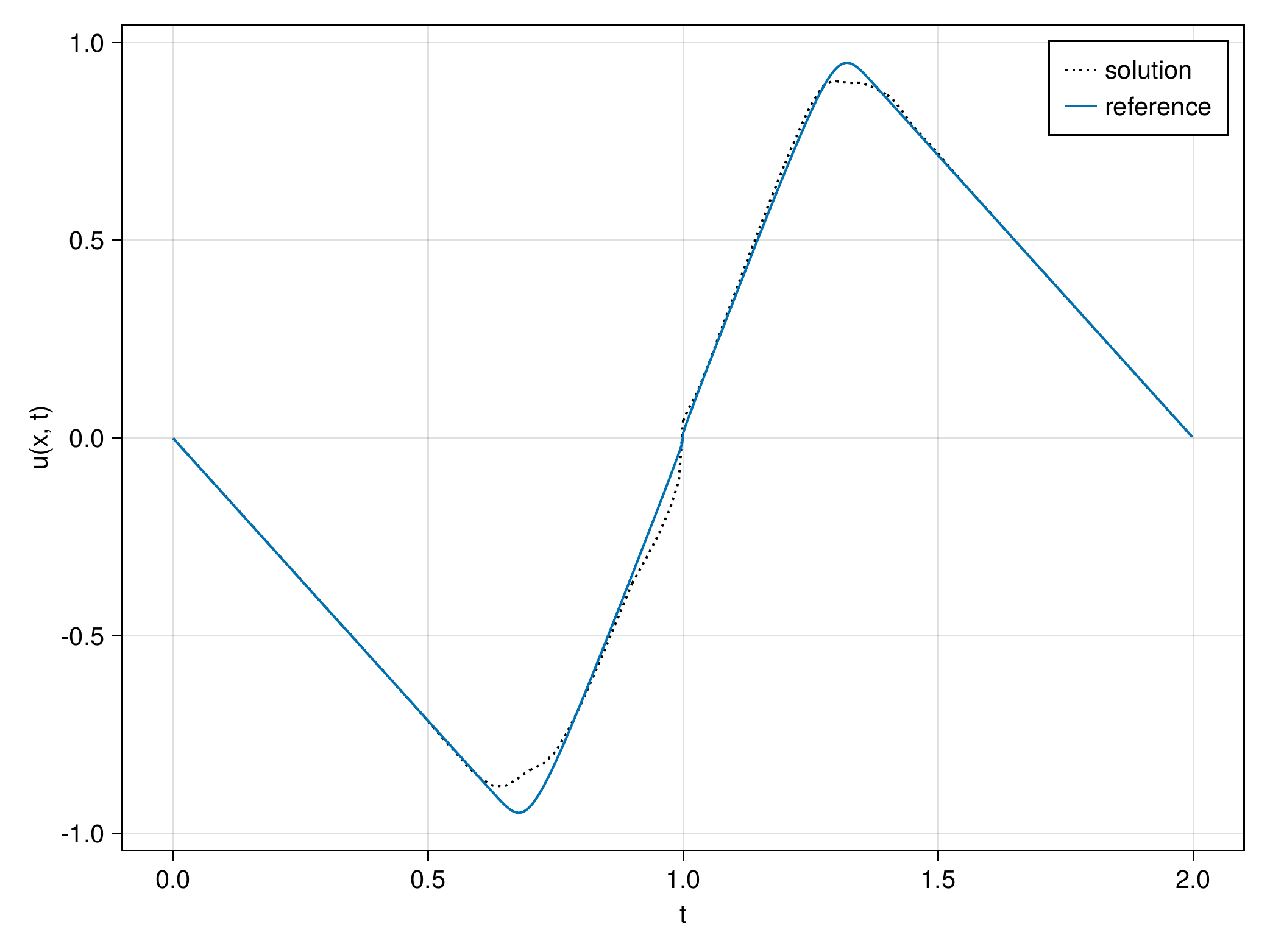}
		\caption{t = 0.3, 20 Elements}
	\end{subfigure}
		\begin{subfigure}{0.49 \textwidth}
		\includegraphics[width=\textwidth]{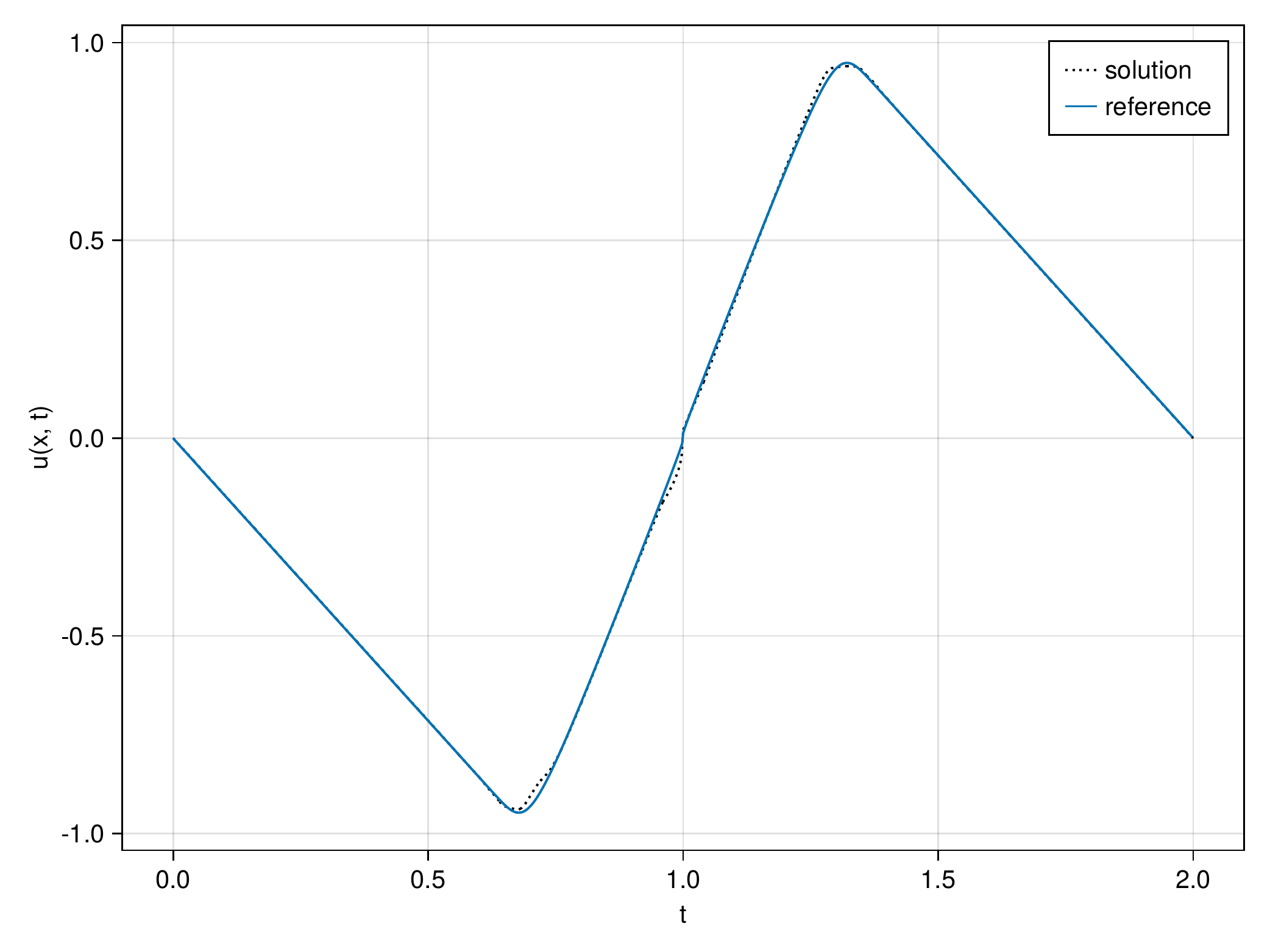}
		\caption{t = 0.3, 50 Elements}
	\end{subfigure}
	\caption{Solution to the second initial condition calculated by the DDG scheme.}
	\label{fig:raretest}
\end{figure}
\begin{figure}
	\begin{subfigure}{0.49 \textwidth}
		\includegraphics[width=\textwidth]{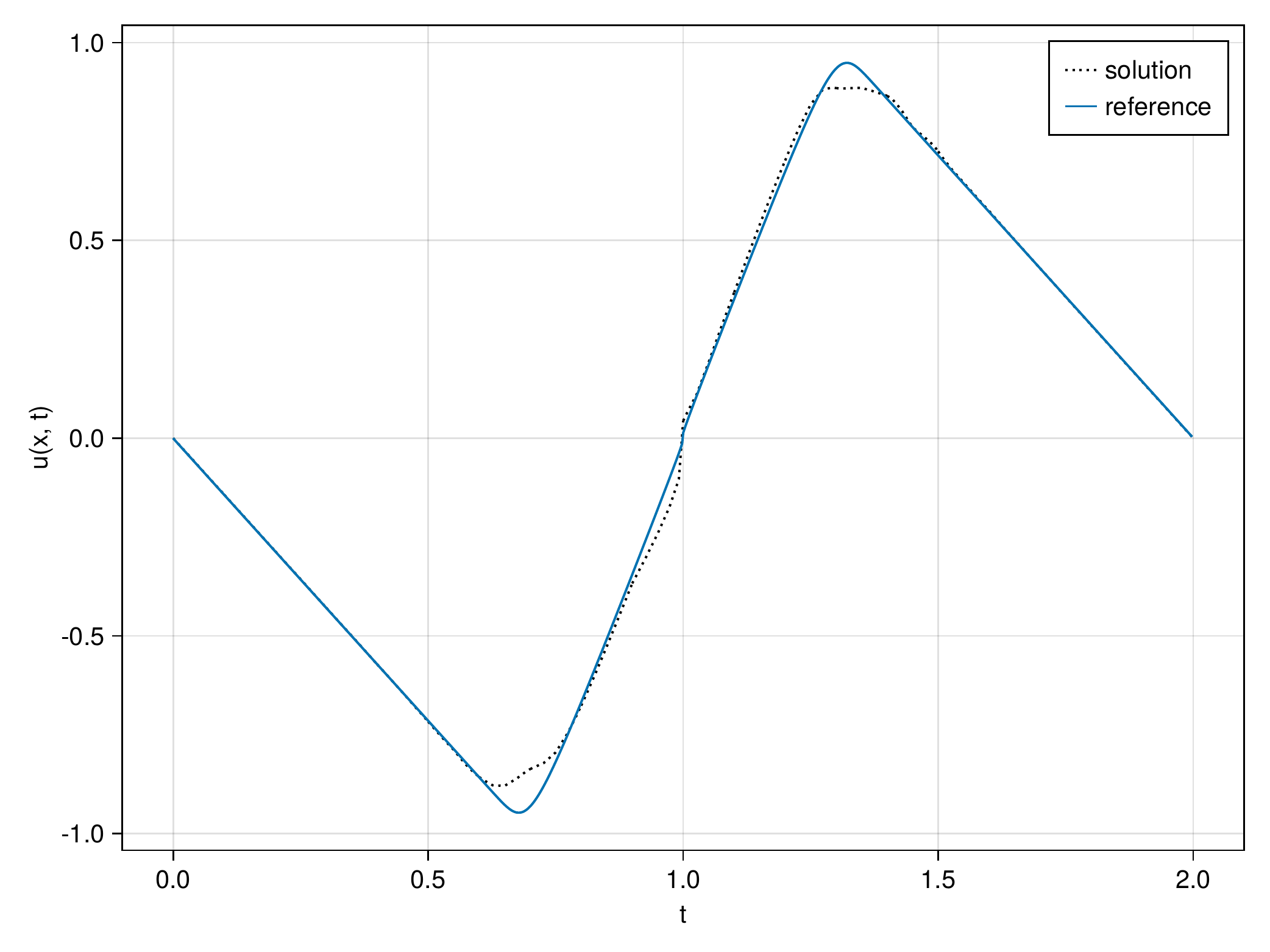}
		\caption{t = 0.3, 20 Elements}
	\end{subfigure}
	\begin{subfigure}{0.49 \textwidth}
		\includegraphics[width=\textwidth]{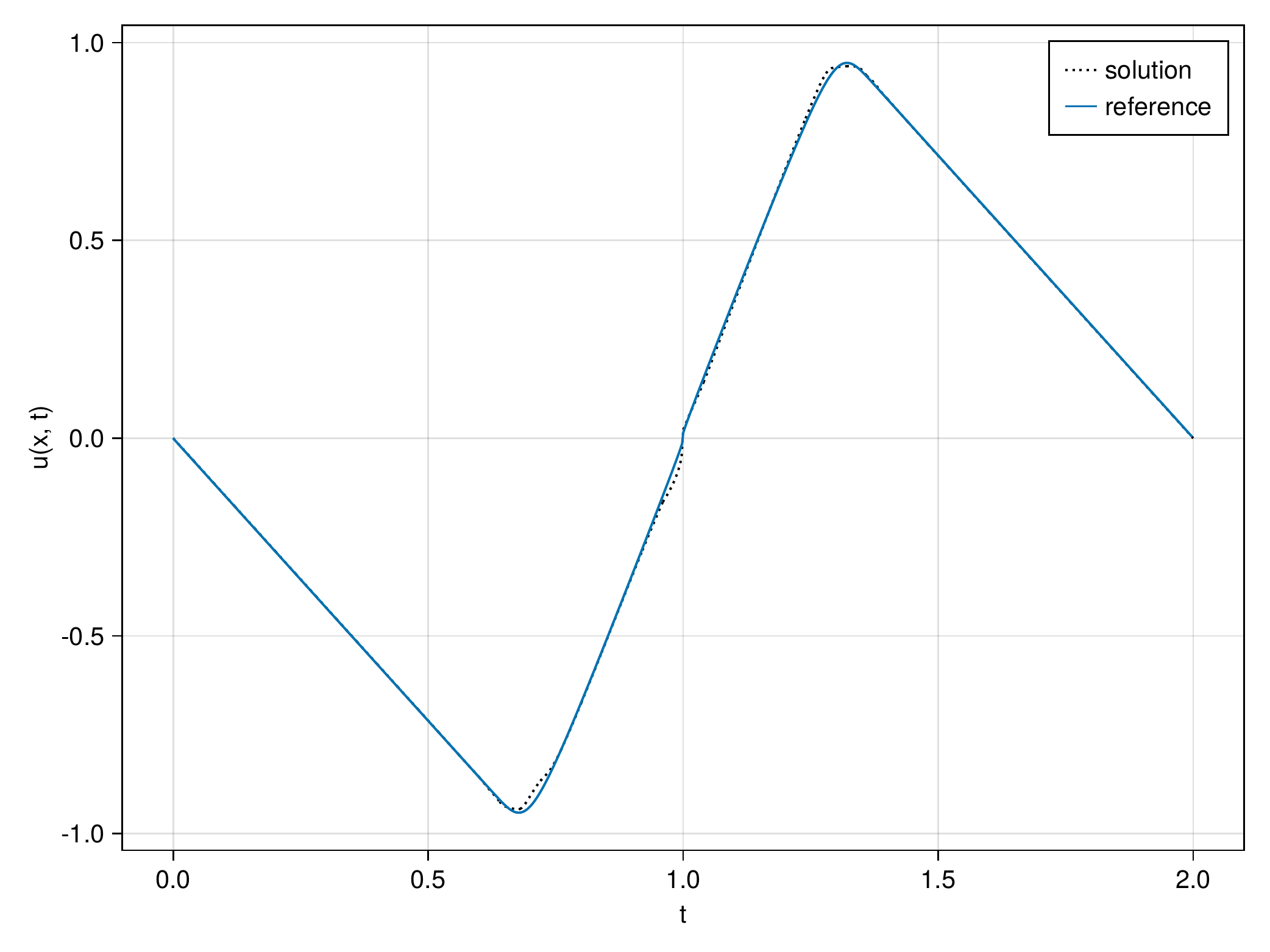}
		\caption{t = 0.3, 50 Elements}
	\end{subfigure}
	\caption{Solution to the second initial condition calculated by the fully discrete scheme DRKDG.}
	\label{fig:draretest}
\end{figure}
First of all, both test cases were successfully run by both schemes up to any simulation time. Moreover, as can be seen from the results in Figure \ref{fig:sintest} and Figure \ref{fig:dst}, their solutions seem to be essentially free of oscillations, as only the polynomials in up to three shocked cells oscillate slightly. This evidence of a robust scheme is only hampered slightly by the tests concerning the sonic point glitch. The cells around the sonic point of the flux clearly show a problematic feature of the base Lax-Friedrichs scheme that can't be corrected using the devised method, consult also Figure \ref{fig:raretest} and Figure \ref{fig:draretest}. Future improvements could therefore be based around modifications of the intercell flux.
\subsection{Numerical test of the Dafermos entropy rate criterion and semidiscrete entropy inequalities}
In \cite{klein2022using} the author tested several schemes for their compatibility with Dafermos' entropy rate criterion. A followup paper  \cite{klein2022Comparison} also tested if a similar family of solvers respects the classical entropy inequality, i.e.~if the schemes are entropy dissipative. We would like to test the new DG scheme presented in this publication for the same two entropy criteria, i.e. Dafermos' entropy rate criterion and classical entropy inequalities. Our first test case with initial condition $u_1(x, 0)$ will be used, which coincides with the test case used in the first publication.
\begin{figure}
	\begin{subfigure}{0.49 \textwidth}
		 \includegraphics[width=\textwidth]{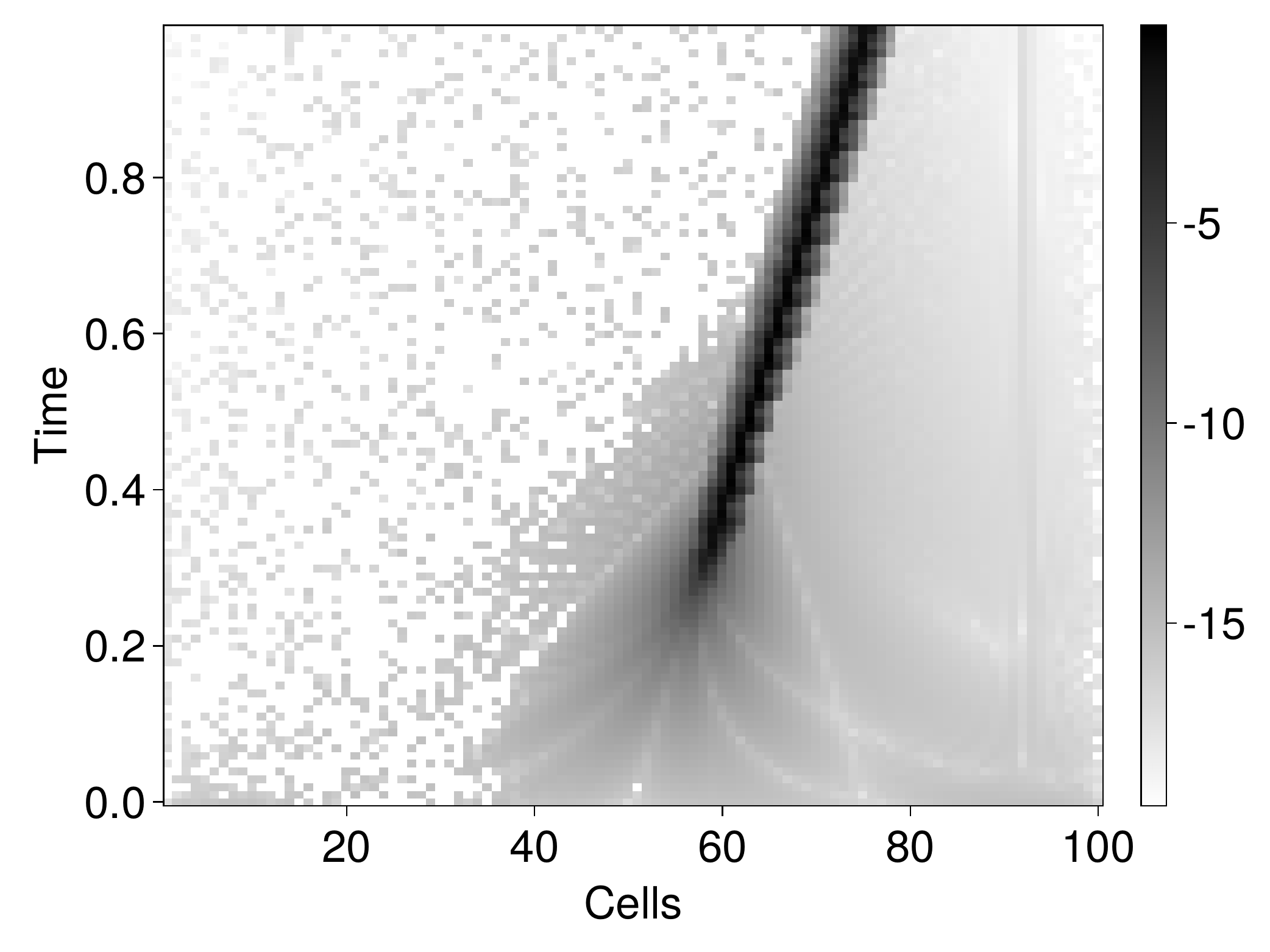}
		 \caption{Logarithm of the negative violation of the entropy equality}
	\end{subfigure}
	\begin{subfigure}{0.49 \textwidth}
		\centering
		\includegraphics[width=\textwidth]{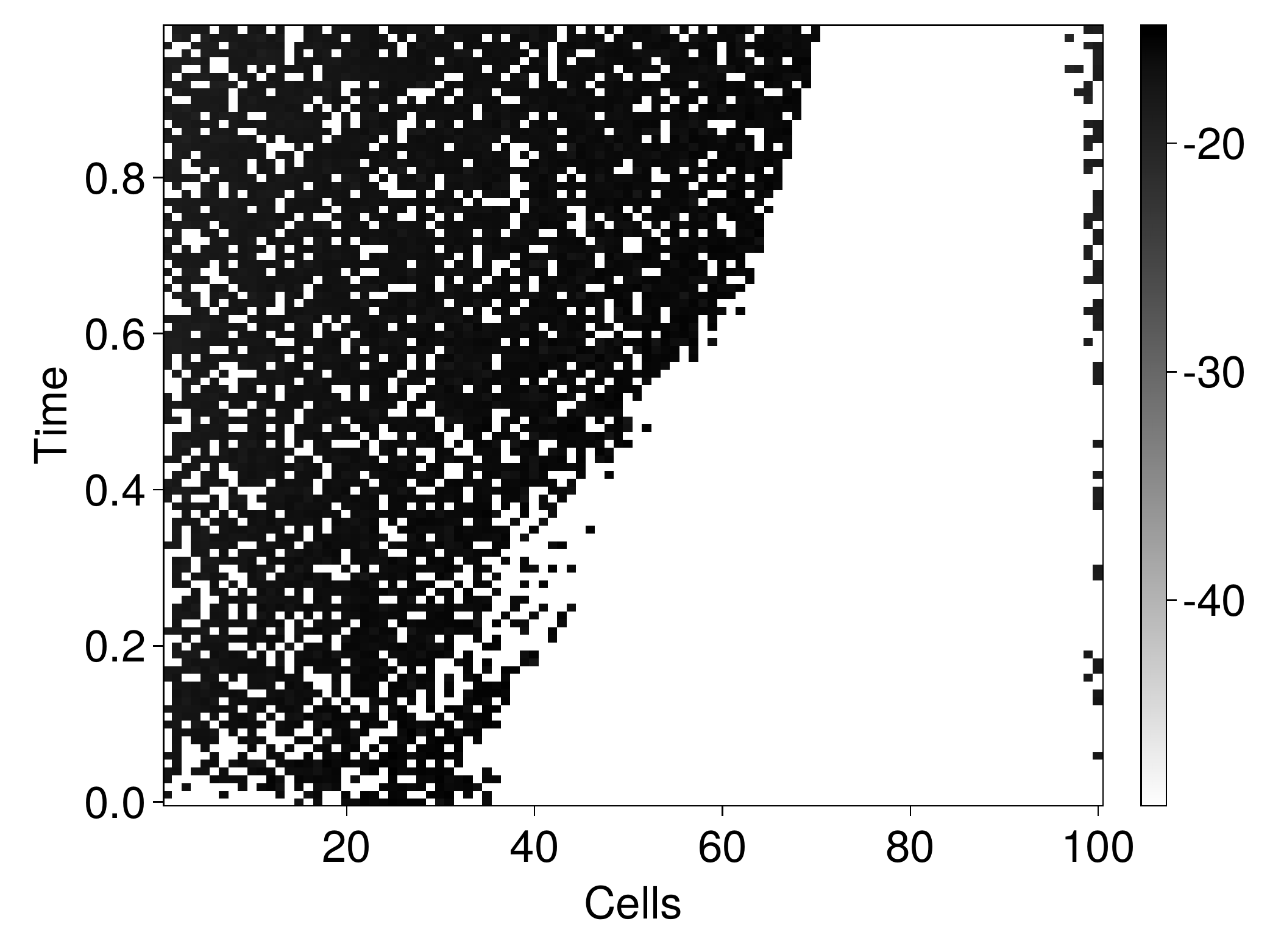}
		\caption{Logarirthm of the negative violation of the entropy equality}
	\end{subfigure}
\caption{Negative and positive violation of the entropy equality for testcase $u_1$. The positive violation is of the same magnitude as the machine precision.}
	\label{fig:Etest}
\end{figure}
\begin{figure}
	\includegraphics[width=0.9\textwidth]{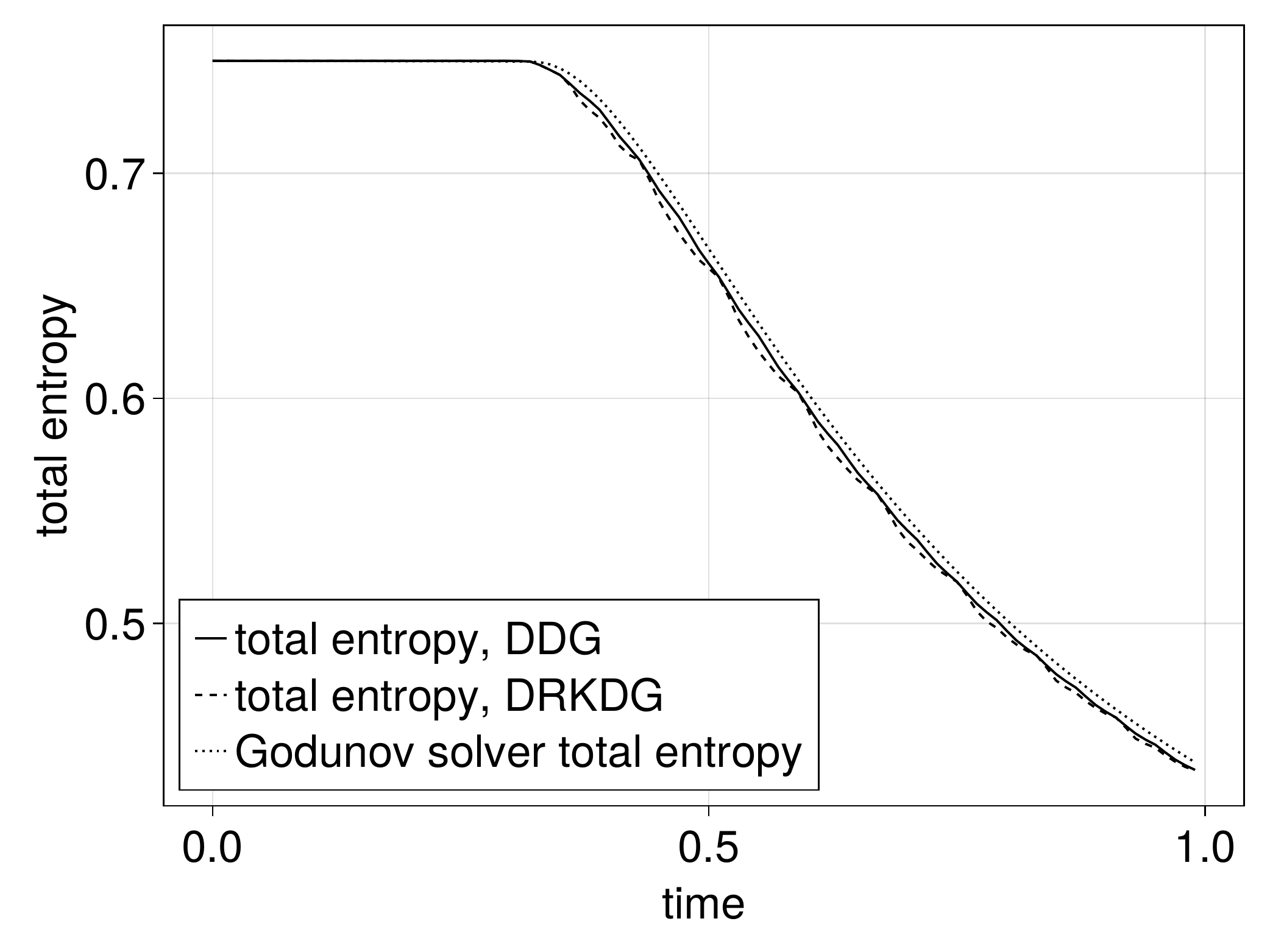}
	\caption{Total entropy of the solution calculated by the DDG and DRKDG schemes in comparison with the total entropy of a solution calculated by a Godunov solver}
	\label{fig:daftest}
\end{figure}
The semi-discrete entropy inequality was calculated for the first test case at every time-step and the decadic logarithm of the positive and negative deviation from the entropy equality was plotted as a heat map, see Figure \ref{fig:Etest}. This was deemed necessary as error bounds stemming from numerical quadrature are used in the implementation, hence the proved entropy inequality does not apply in a strict sense but under the assumption of vanishing quadrature errors. Still, the hope is that the deviations from this entropy equality are small. And interestingly this seems to be true. A positive violation is only appearing with a magnitude of $10^{-16}$, around the machine precision used for the calculation and only in extremely smooth areas of the solution - indicating that these smooth areas indeed lead to an entropy equality, as dictated by theory. Moreover, we would like to find evidence of a numerical solution aligned with Dafermos entropy rate criterion. Therefore the same method as in \cite{klein2022using} was used. The total entropies of both solvers can be seen in Figure \ref{fig:daftest}. The solution of an extremely fine grained Godunov scheme with $N = 10000$ cells was added as a reference solution, to which the entropy of the numerical solutions by the modified DG schemes can be compared to. Both curves align nearly perfectly, with one exception. The entropy of the DG method starts to drop at nearly the same rate as the one of the Godunov method, but at a slightly earlier time. While the rates, apart from small oscillations in the DG method, are comparable, the total entropy of the DG method is always lower than the one of the Godunov method. Against the first instinct, this is not an unpleasant result. The reason being that the DG method is based around a piecewise polynomial approximation space. Therefore the method starts every timestep from an approximate solution and this solution seems to be unable to resolve the sharp discontinuity and instead smears it over up to $3$ cells. This, together with conservation, seems to force a lower total entropy for the approximate solution. One can therefore conclude that our plan to enforce Dafermos' entropy rate criterion in DG schemes was successful.

\subsection{Numerical convergence analysis}
Our deviations from the basic DG scheme are only limited by an assumed error estimate. Future work could therefore be concentrated on deriving convergence rate estimates for smooth solutions. We will collect some numerical evidence that such convergence rate estimates are possible. This will be done by an experimental convergence analysis. As we will see, our schemes retain the high order of accuracy of the base schemes. Once more Burgers' equation was solved with the initial condition
\[
	u(x, 0) = 1 + \frac{\sin(\pi x)}{50}
\]
and a periodic domain \( \Omega = [0, 2) \). The solution of this problem is smooth on $t = [0, 8]$ and we therefore look at the error at $t_{\mathrm{end}} = 8$ between the solution calculated with $N \in \{10, 15, 20, 25, 30\}$ cells and internal polynomials of order $p=6$. The relatively small maximum amount of cells was chosen as it was feared that otherwise the floating point accuracy could interfere with the calculation of the error estimates. Time integration was carried out once more using the SSPRK33 method. The time step was not selected to keep $\lambda = \frac{\Delta t}{\Delta x}$ constant, but as $\lambda (\Delta x)^2= \Delta t$ as this allow us to observe up to sixth order convergence rate without being limited by the accuracy of the time integration method. The method displays sixth order of accuracy in our tests 
\begin{figure}
	\begin{subfigure}{0.5\textwidth}
	\includegraphics[width=\textwidth]{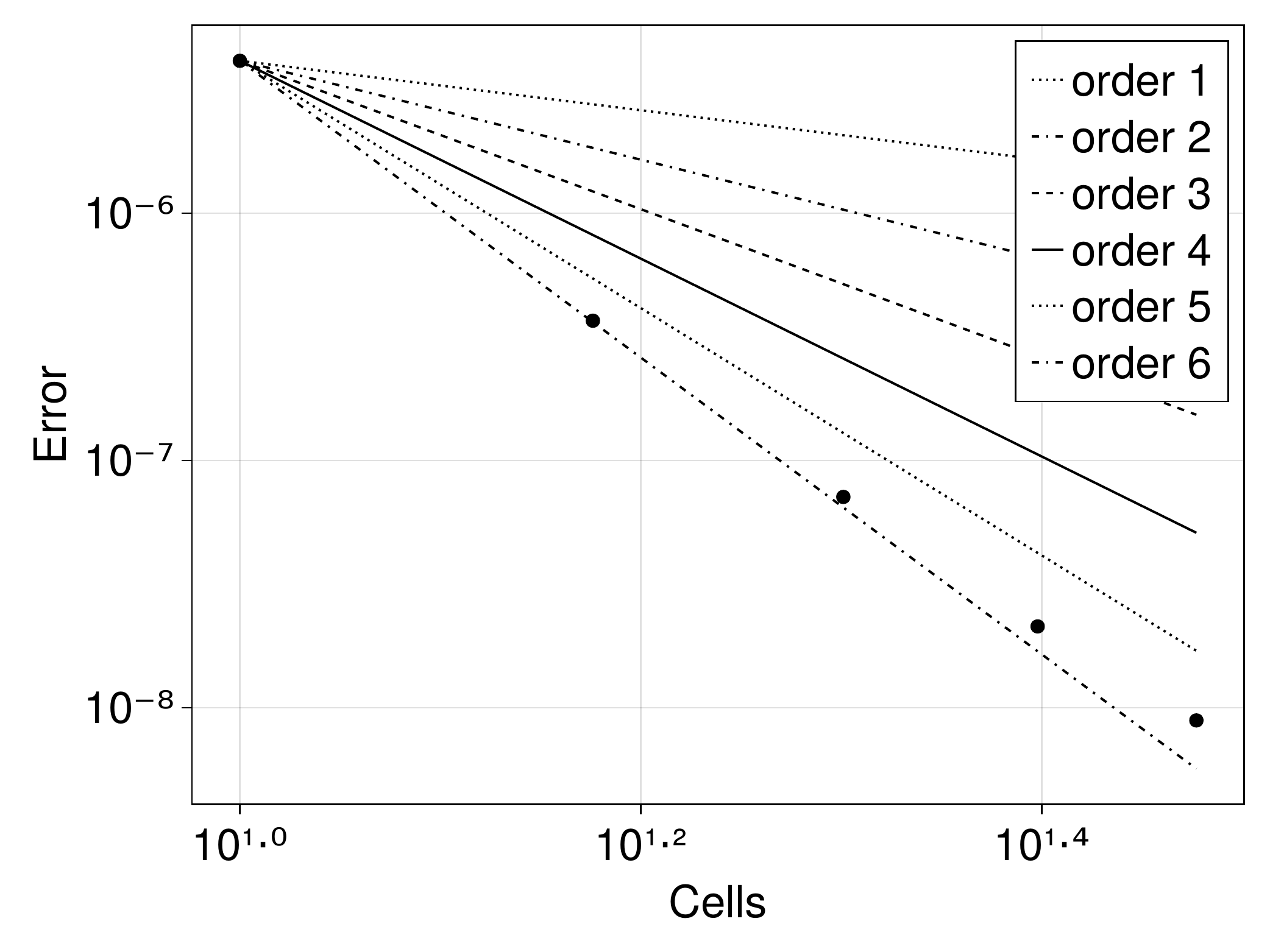}
	\caption{Convergence analysis for polynomial order $6$ in the 1-Norm.}
\end{subfigure}
	\begin{subfigure}{0.5\textwidth}
		\includegraphics[width=\textwidth]{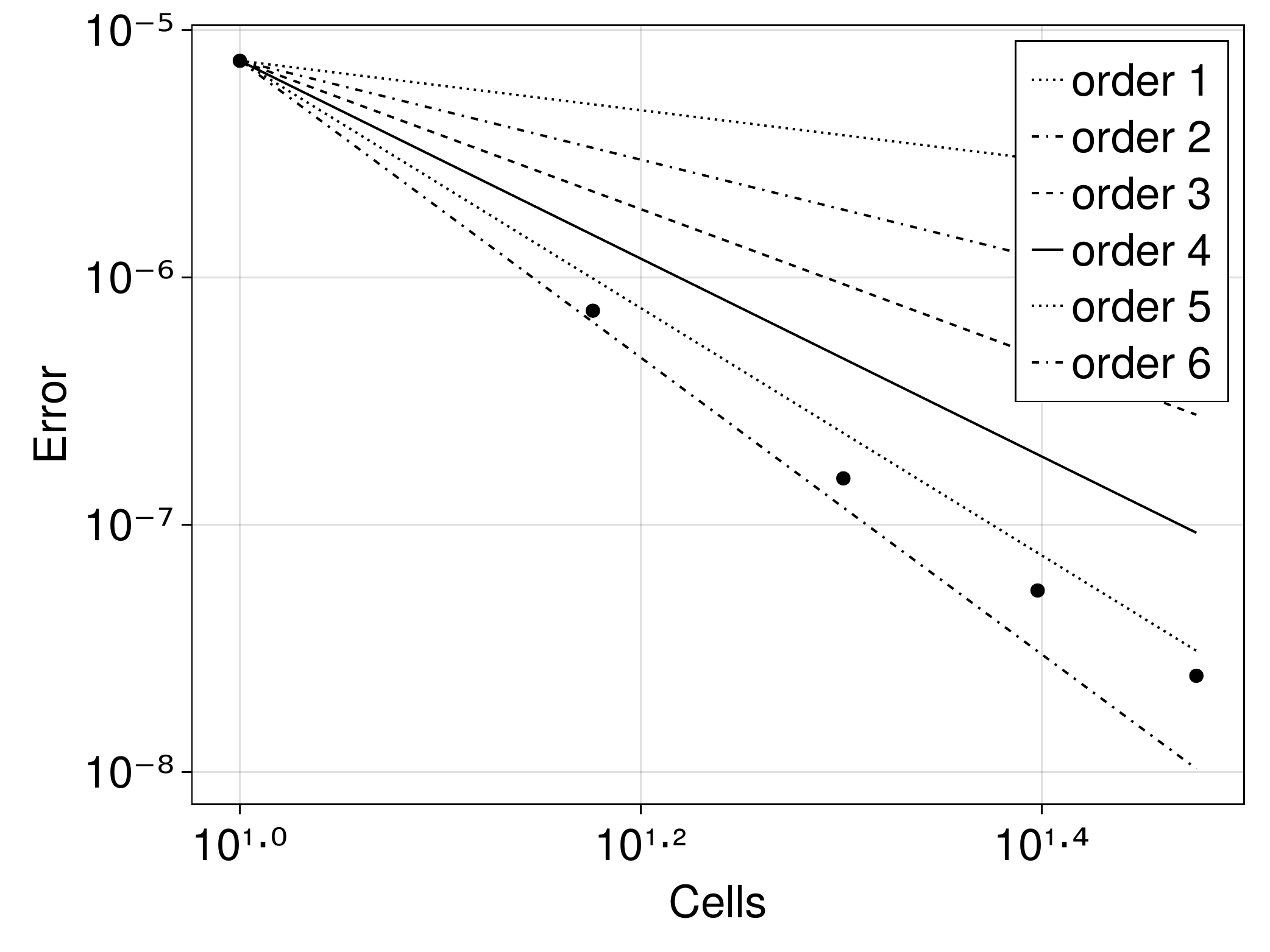}
		\caption{Convergence analysis for polynomial order $6$ in the 2-Norm.}
	\end{subfigure}
		\caption{Convergence analysis for the semi-discrete scheme DDG.}
	\label{fig:CASD}
\end{figure}
\begin{figure}
	\begin{subfigure}{0.5\textwidth}
		\includegraphics[width=\textwidth]{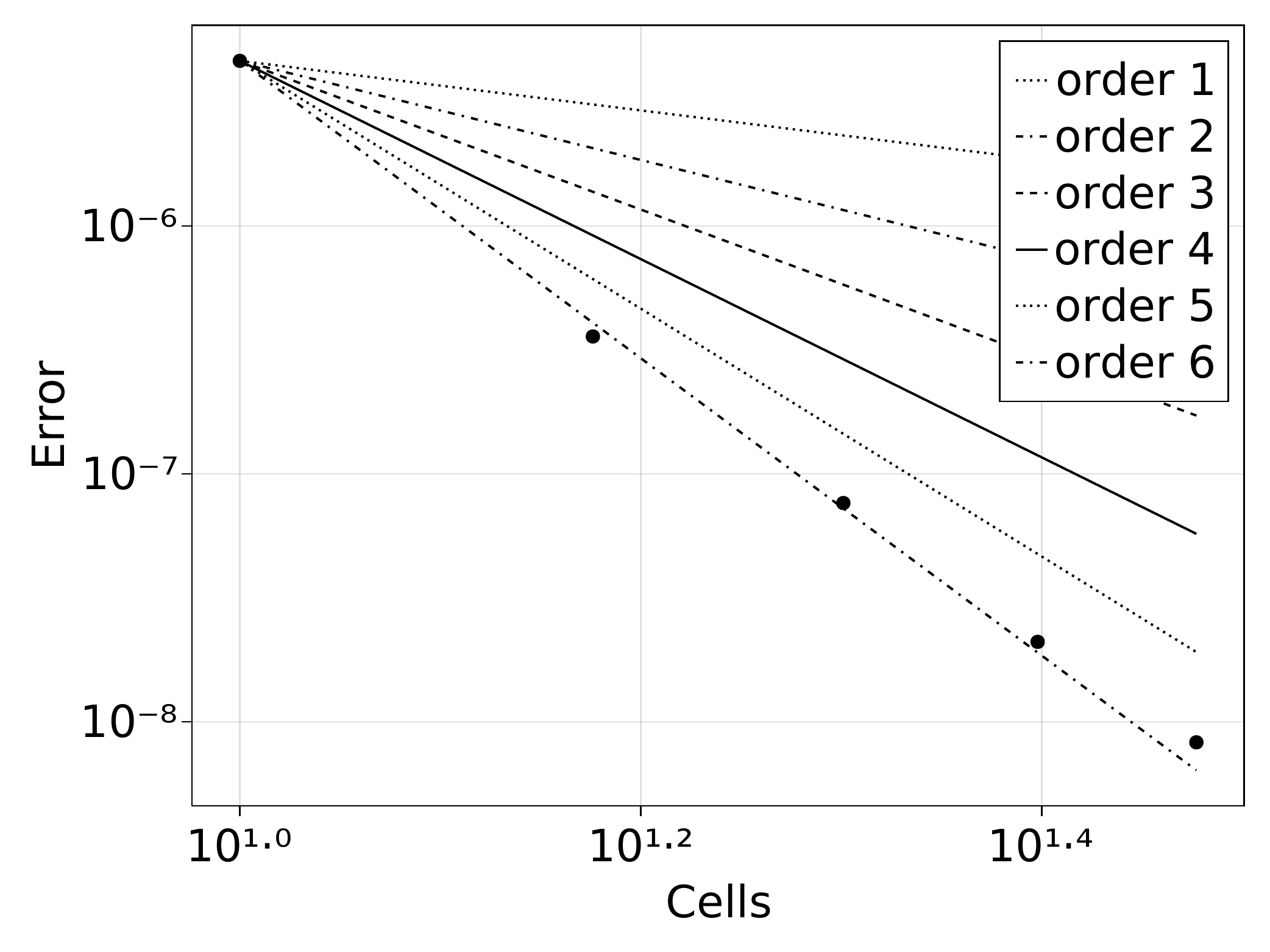}
		\caption{Convergence analysis for polynomial order $6$ in the 1-Norm for the discrete scheme.}
	\end{subfigure}
	\begin{subfigure}{0.5\textwidth}
		\includegraphics[width=\textwidth]{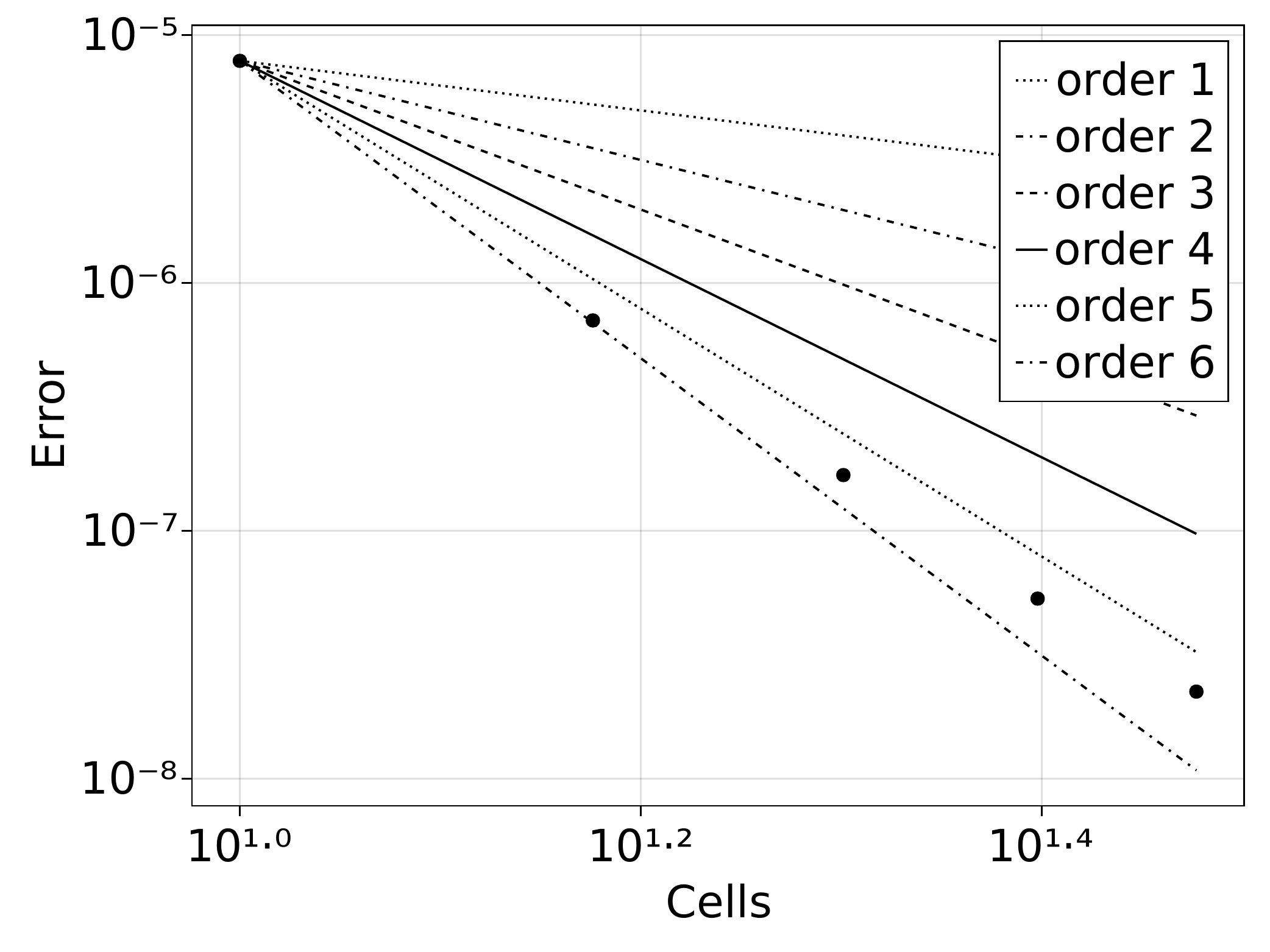}
		\caption{Convergence analysis for polynomial order $6$ in the 2-Norm for the discrete scheme.}
	\end{subfigure}
	\caption{Convergence Analysis for the fully discrete scheme DRKDG.}
\label{fig:CAFD}	
\end{figure}
in Figures \ref{fig:CASD} and \ref{fig:CAFD}. While this is a promising result some caution is advised. The author witnessed problems concerning the high order convergence of lower order versions of the scheme. These problems, for example present for order $3$, could have several reasons. One of them being that the method was implemented without any modal filtering present. Higher polynomial orders could have lead to smaller aliasing errors and therefore mitigated the problem in the tests. Future tests will include modal filtering as an additional building block after the theoretical compatibility of modal filtering and the DG entropy descent was explored. See \cite{GOS2018Modal, RGOS2018Stab} concerning modal filtering.

\subsection{Experimental analysis of the timestep restriction}
Sadly, the vanilla DG method suffers from low time step restrictions, i.e.~ the CFL number, the upper bound on the grid constant $\lambda = \frac{\Delta t}{\Delta x}$ with respect to the highest signal velocity, drops drastically with the used order \cite{Cockburn1989DGI, CockburnShu1989DGI, HW2008DG}. Several methods have been devised to counteract this problem and since the limiting base FV scheme of our discretization does not have such a low time step restriction we will test which effect our modifications have on the allowed time step. We use the first initial condition given, that was used as a test for discontinuous solutions. Note that the standard DG scheme blows up at $t \approx 0.3$ when a shock forms independently of the time step. We test the scheme with different CFL numbers and orders. We can therefore identify, if the maximum allowed time step scales linearly with the grid size and if, and how their ratio depends on the order of the used cells. 
\begin{figure}
	\begin{subfigure}{0.49\textwidth}
		\includegraphics[width=\textwidth]{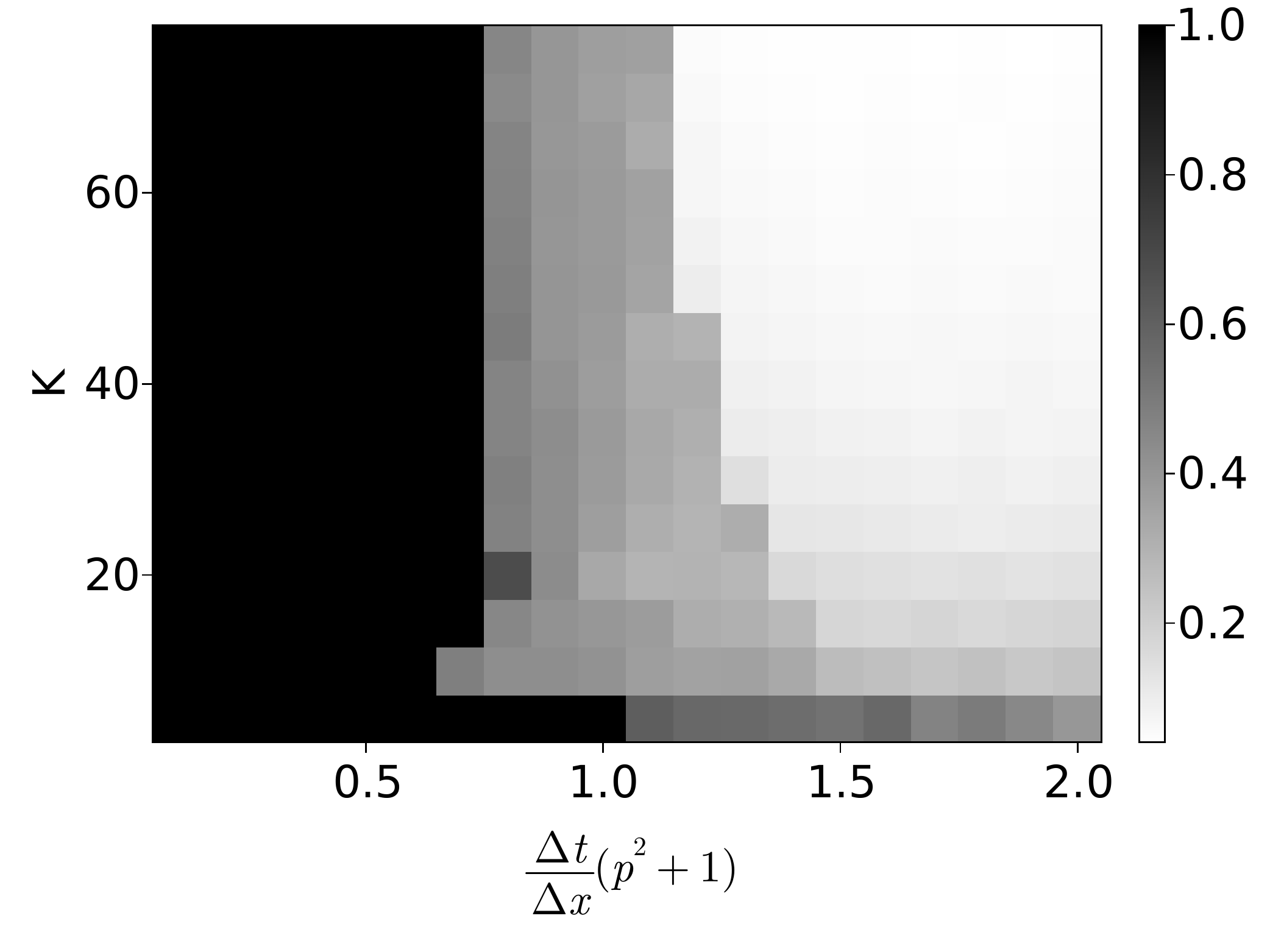}
		\caption{p=3}
	\end{subfigure}
	\begin{subfigure}{0.49\textwidth}
		\includegraphics[width=\textwidth]{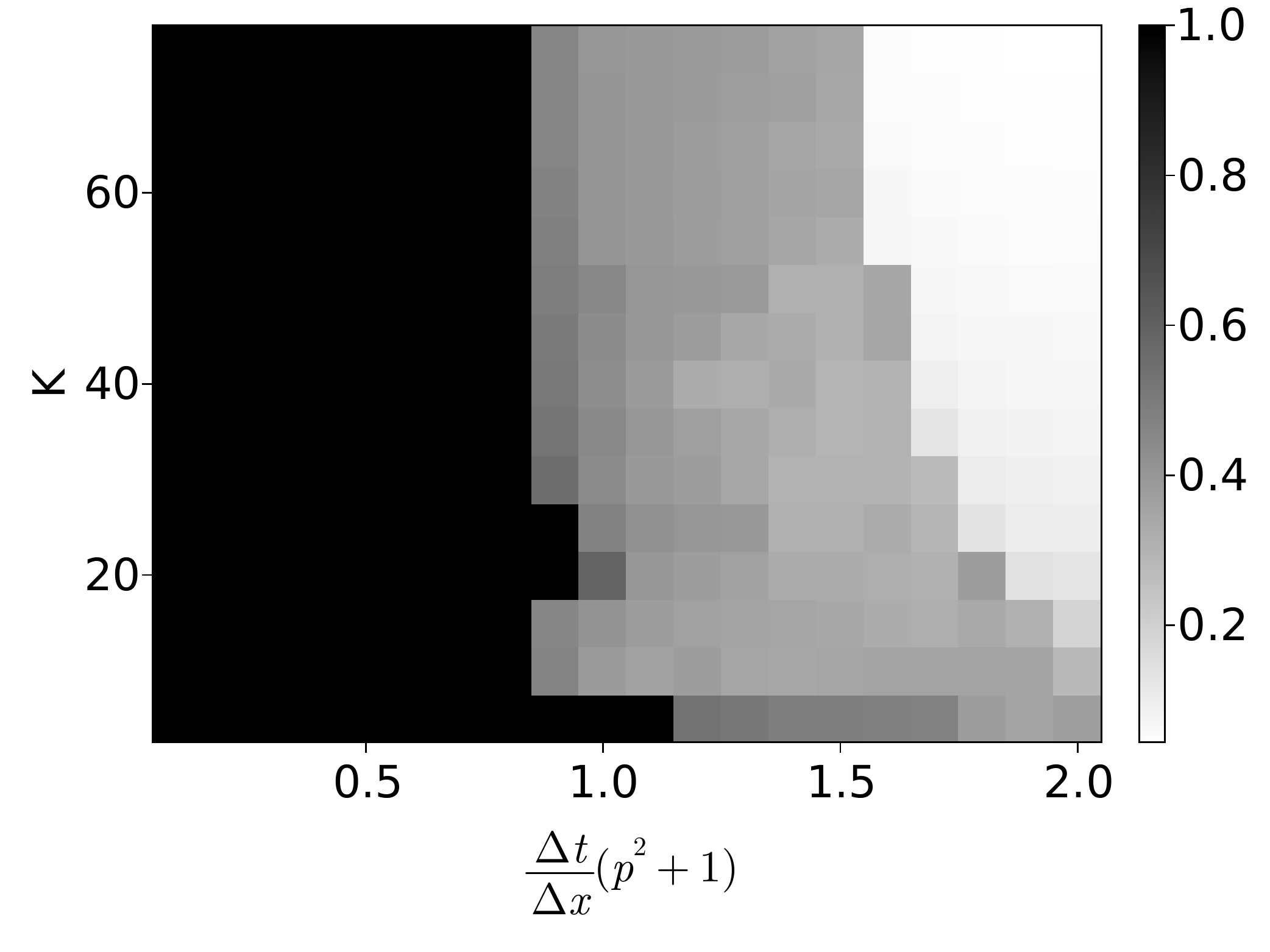}
		\caption{p=4}
	\end{subfigure}
	\begin{subfigure}{0.49\textwidth}
		\includegraphics[width=\textwidth]{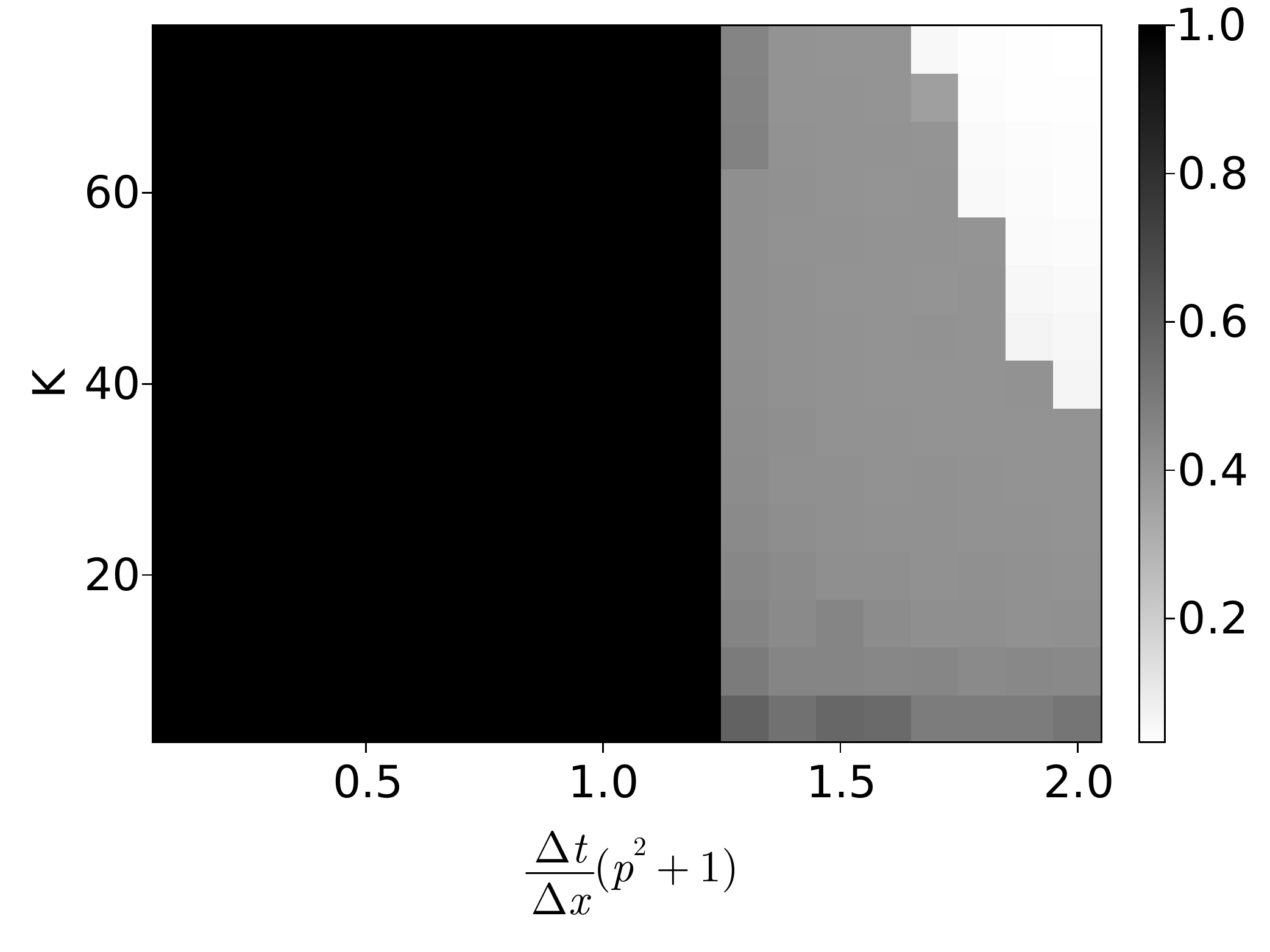}
		\caption{p=5}
	\end{subfigure}
	\begin{subfigure}{0.49\textwidth}
		\includegraphics[width=\textwidth]{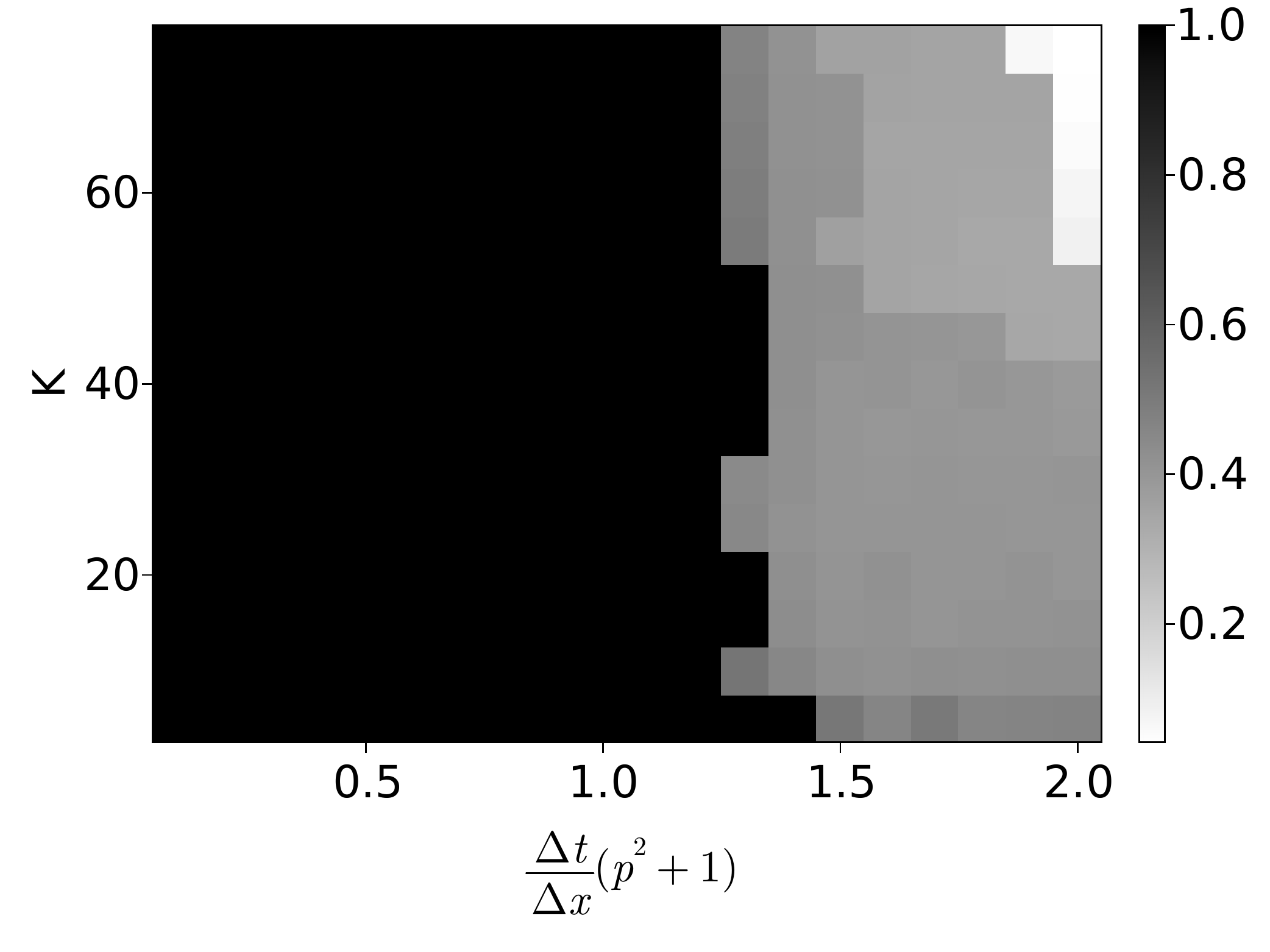}
		\caption{p=6}
	\end{subfigure}
\caption{Heatmap of the achieved maximum simulation time before the solution blows up for test case 1 in relation to the used CFL number and the number of cells. The semidiscrete scheme DDG was used.}
\label{fig:blowup}
\end{figure}
\begin{figure}
	\begin{subfigure}{0.49\textwidth}
		\includegraphics[width=\textwidth]{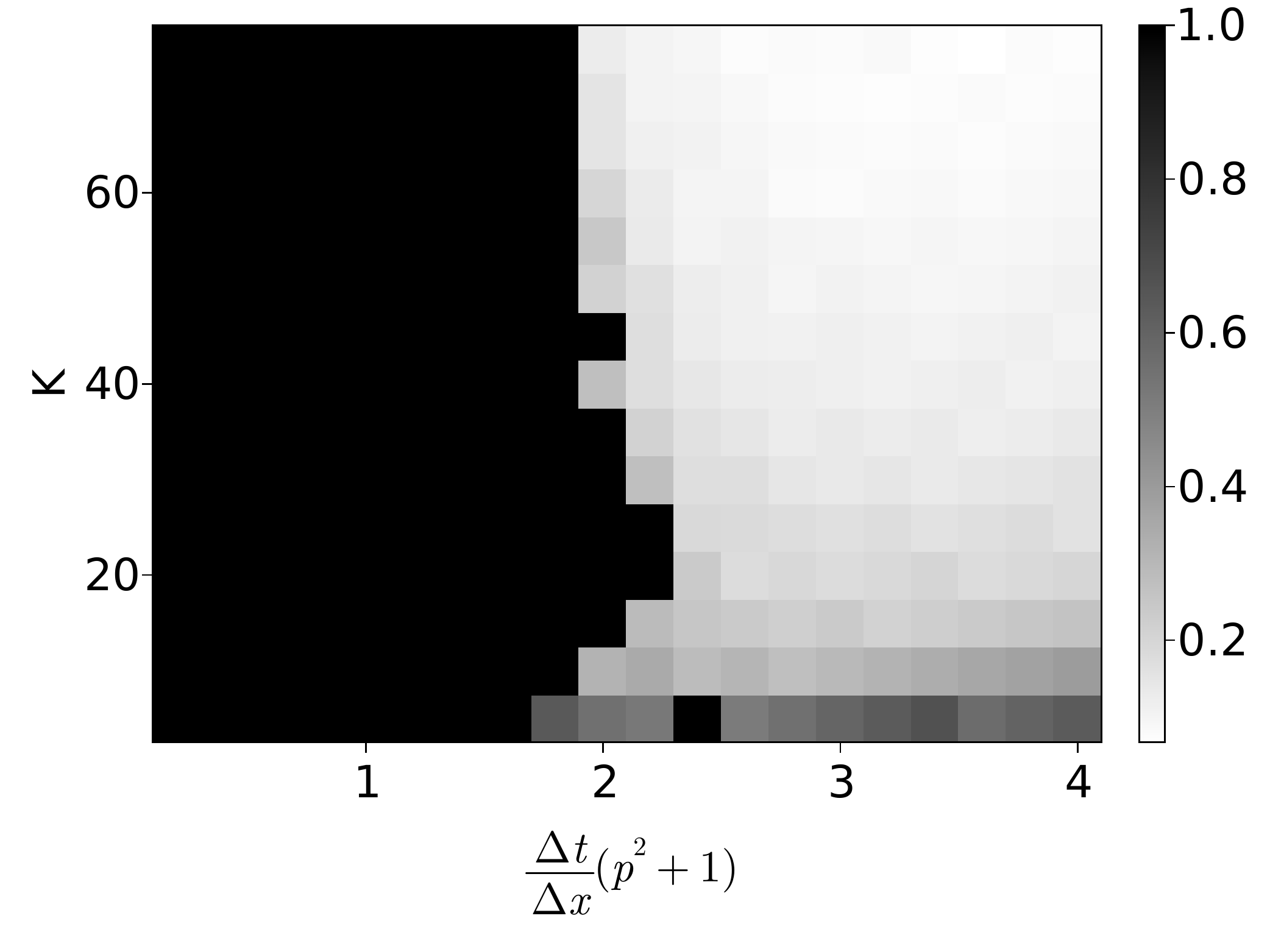}
		\caption{p=3}
	\end{subfigure}
	\begin{subfigure}{0.49\textwidth}
		\includegraphics[width=\textwidth]{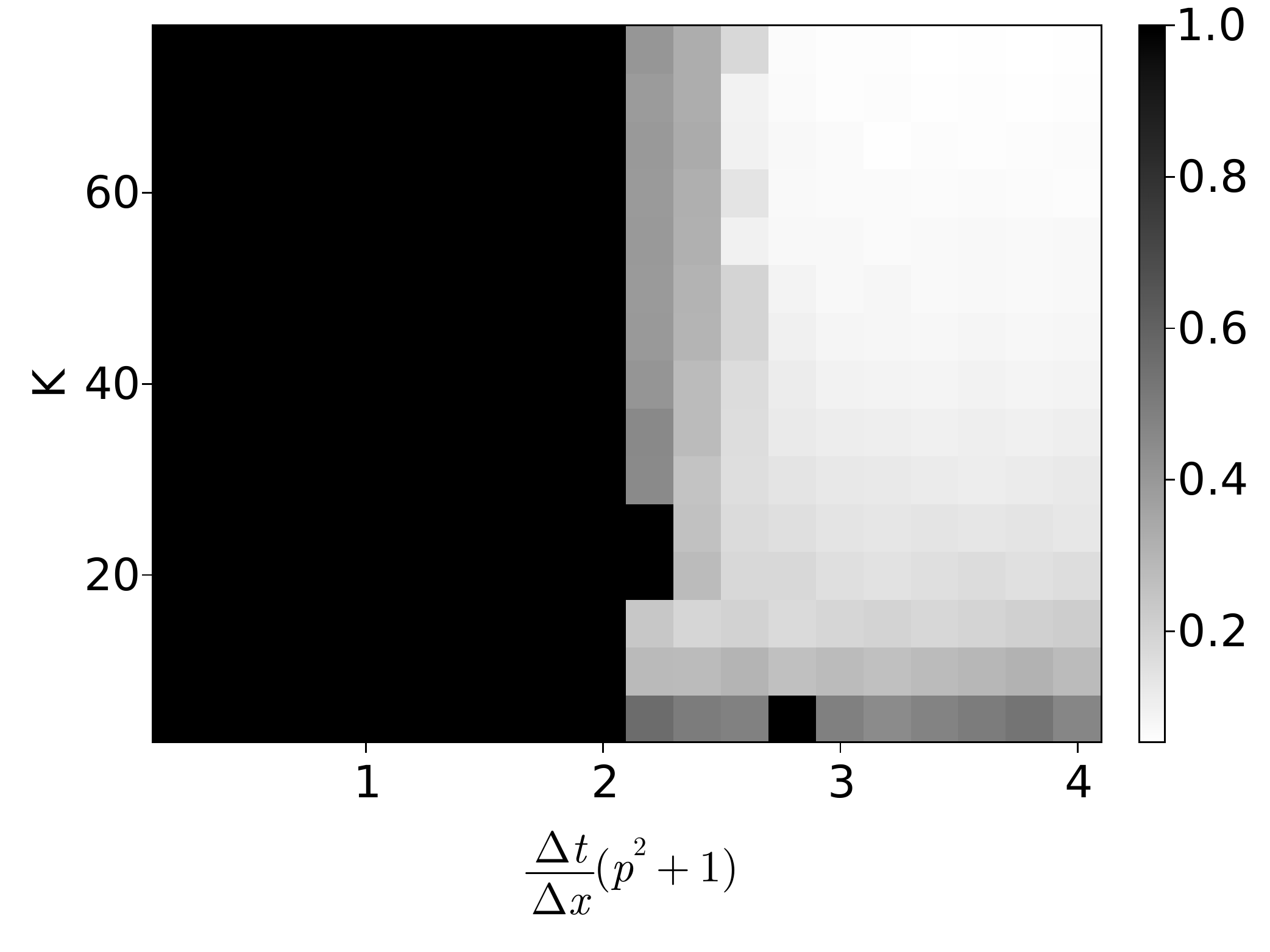}
		\caption{p=4}
	\end{subfigure}
	\begin{subfigure}{0.49\textwidth}
		\includegraphics[width=\textwidth]{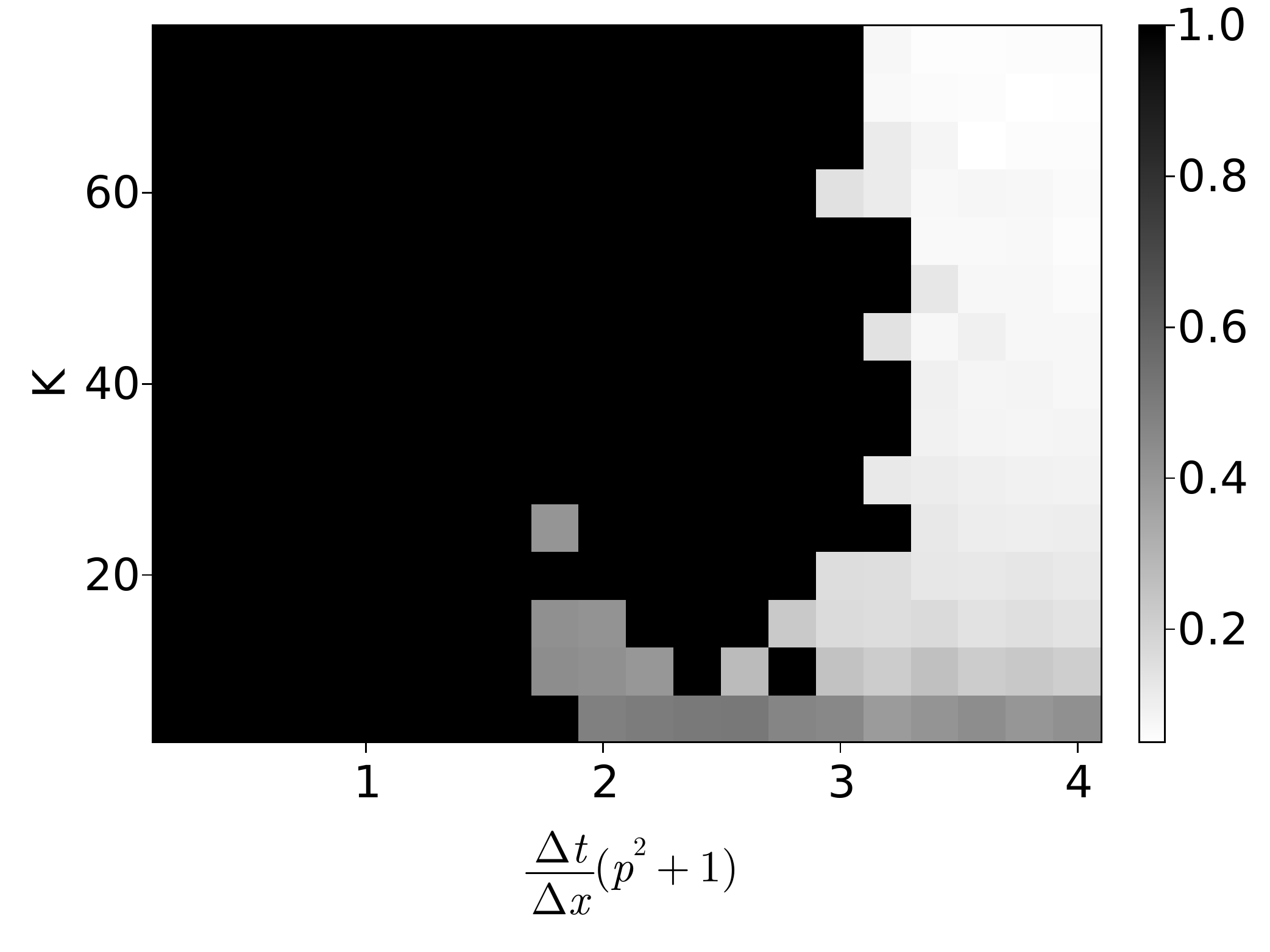}
		\caption{p=5}
	\end{subfigure}
	\begin{subfigure}{0.49\textwidth}
		\includegraphics[width=\textwidth]{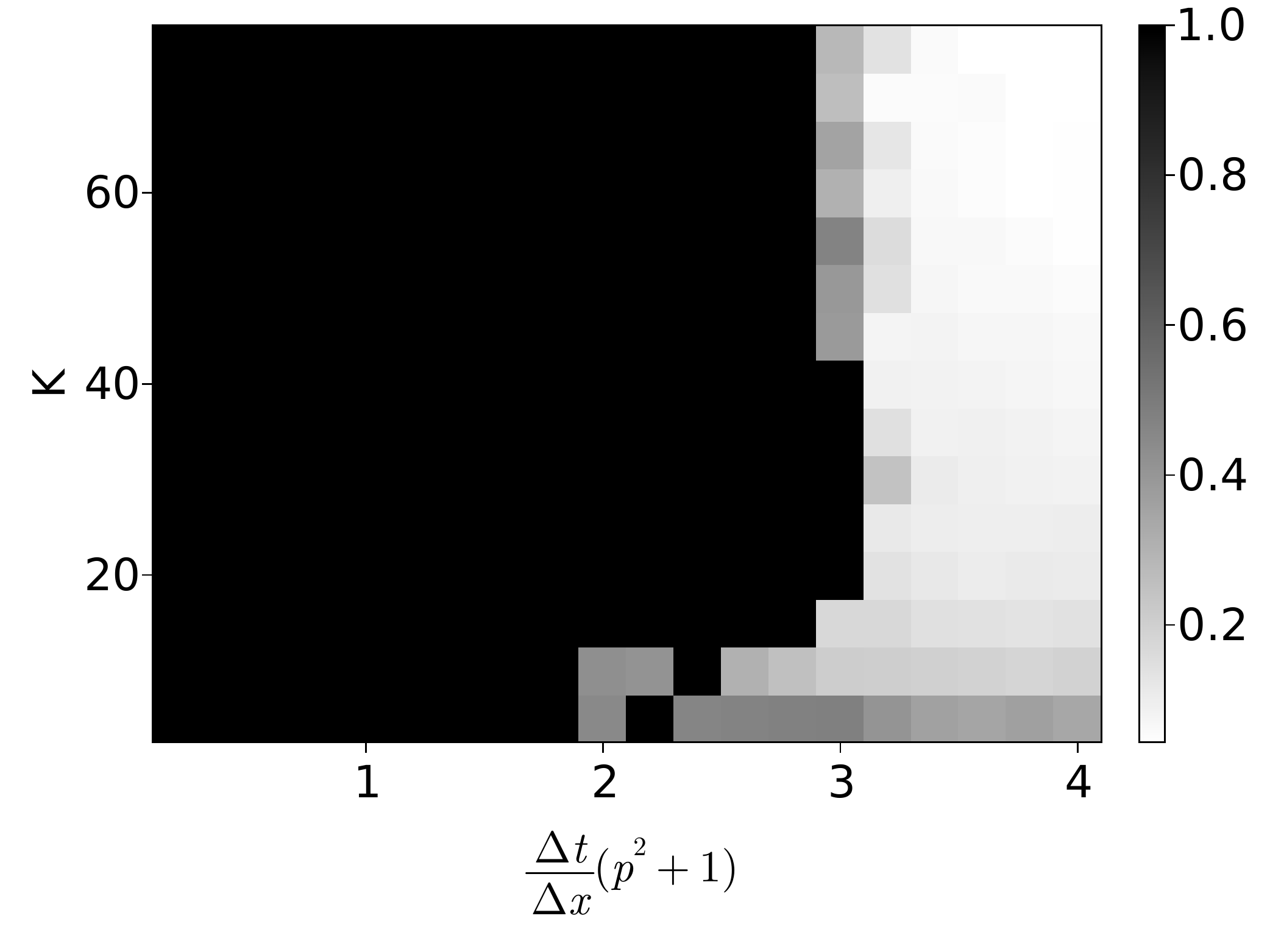}
		\caption{p=6}
	\end{subfigure}
	\caption{Heatmap of the achieved maximum simulation time before the solution blows up for test case 1 in relation to the used CFL number and the number of cells for the discrete scheme DRKDG}
	\label{fig:dblowup}
\end{figure}
 
The results for our modified DG schemes are shown in Figure \ref{fig:blowup} and Figure \ref{fig:dblowup}. The shown heatmaps correspond to the simulation time at which a blowup occurred. The black areas correspond to the maximum simulation time of $t = 1$. There, no blowup occurred up to this time, where the simulation was stopped. Interestingly, the DRKDG scheme is able to use significantly bigger time steps than the DDG scheme, when the Grid is fine enough. Gains of a factor of $4$ for polynomials of degrees 3 to 6 can be seen in Figure \ref{fig:dblowup} compared to  Figure \ref{fig:blowup}. A problem in this regard concerns the accuracy of solutions calculated in this way. The simulations do not crash, but are pulled down to first order when used at to big time steps. This is especially dangerous if these schemes are used with the rule of thump to run them with time steps only barely stable.
	\section{Conclusion} \label{sec:concl}
	In this publication, the author used error estimates between the exact entropy variables and the approximate entropy variables, error estimates between a fine sub cell scheme and the approximate solution, and error estimates between a projection of the exact solution to the DG approximation space to control the entropy in DG methods. It was thereby possible to bound the derivative of the total entropy with the derivative of the total entropy of the limit of the sub cell scheme. This serves as a numerical treatment of the Dafermos entropy criterion in its original form with the sub cell solution as a reference weak solution to which the modified DG solution is compared concerning its total entropy dissipation. If one conjects that this solution already satisfies the Dafermos criterion, it follows further that the approximate solution of the DG method satisfies this criterion. Further, the author conjectured in a previous work that a numerical Dafermos entropy rate criterion should be centered around allowing only a minimization of the entropy rate if the additional residual produced by this action is small. This is especially important as the method would otherwise be able to dissipate entropy to fast by enlarging approximation errors. Because the scheme was designated as a correction to the classical DG method, where the correction is the biggest possible entropy dissipation in the allowed error margin, could the method be also interpreted as correctly implementing this criterion.
As a side note a classical entropy inequality is also satisfied if one conjectures the correctness of the error estimators and scales the maximum descent direction using the formulae derived by the author. The same technique was also used to construct a fully discrete variant of the scheme. Apart from being able to calculate smooth solutions with high order of accuracy, as the standard DG method, shocks are also handled well without any additional stabilization. Especially oscillations in shocked cells are nearly invisible and do not produce oscillations in other cells. Future work will focus on the application to multiple dimensions and systems of conservation laws. The author is sure that the derived error estimates can be also useful for other techniques in use with DG methods. These could be for example positivity preserving correction terms or artificial viscosity.
	\section{Acknowledgements}
	The author would like to thank Thomas Sonar for discussions concerning error estimates for Finite Volume methods. The author was funded under project SO 363/14-1 by the Deutsche Forschungsgemeinschaft (DFG).

	\section{Bibliography}
	\bibliographystyle{plainnat}
	\bibliography{lit}
\end{document}